\documentclass{article}
\usepackage{latexsym,a4}
\makeatletter
\@addtoreset{equation}{subsection}

\makeatother

\newcommand{\separate}{\medskip\noindent}

\def\newsection{ \separate
   \refstepcounter{subsection} 
   {\large\bf \thesubsection\kern.3em}
}
\def\mytheorem#1{
   \separate{\large\bf Theorem#1:\kern.3em}} 

\def\mylemma#1{
   \separate{\large\bf Lemma#1:\kern.3em}} 

\def\mycorollary#1{
   \separate{\large\bf Corollary#1:\kern.3em}} 

\def\myprop#1{
   \separate{\large\bf Proposition#1:\kern.3em}} 

\def\myremark#1{
   \separate{\large\bf Remark#1:\kern.3em}} 

\def\mydefinition#1{
   \separate{\large\bf Definition#1:\kern.3em}} 

\def\Proof{\separate\underline{Proof:}\kern1em}

\newcommand{\GANZ}{{\sf Z\hspace*{-0.4em}Z}}
\newcommand{\REELL}{{\setlength{\unitlength}{1em}
                     \begin{picture}(0.75,1)
                     \put(0,0){\line(0,1){0.69}}
                     \put(0,0){\sf R}
                     \end{picture}
                   }}
\newcommand{\KOMPLEX}{{\setlength{\unitlength}{1em}
                     \begin{picture}(0.7,1)
                     \put(0.34,0){\line(0,1){0.65}}
                     \put(0,0){\sf C}
                     \end{picture}
                   }}

\newcommand{\FGANZ}{\mbox{\tiny{\rm Z\hspace*{-0.45em}Z}}}

\newcommand{\FKOMPLEX}{\mbox{\tiny{\setlength{\unitlength}{1em}
                               \begin{picture}(0.6,0.5)
                               \put(0.34,0){\line(0,1){0.47}}
                               \put(0,0){\rm C}
                               \end{picture}
                              }}}
\newcommand{\FREELL}{\mbox{\tiny{\setlength{\unitlength}{1em}
                     \begin{picture}(0.6,0.5)
                     \put(0,0){\line(0,1){0.48}}
                     \put(0,0){\rm R}
                     \end{picture}
                   }}}
\newcommand{\Brr}{{\mathchoice{\REELL}{\REELL}{\!\FREELL}{\!\FREELL}}}
\newcommand{\Bcc}{{\mathchoice{\KOMPLEX}{\KOMPLEX}{\!\FKOMPLEX}{
\!\FKOMPLEX}}}

\newcommand{\Bii}{{\mathchoice{\GANZ}{\GANZ}{\FGANZ}{\FGANZ}}}
\newcommand{\CPE}{{\mathchoice{\Bcc{\rm P}_1}{\Bcc{\rm P}_1}{\Bcc\!\!
\mbox{\rm\tiny P}_1}{\Bcc\!\!\mbox{\rm\tiny P}_1}}}

\def\QED{\hfill$\Box$}
\def\inv{^{-1}}
\def\BEA{\begin{eqnarray}}
\def\EEA{\end{eqnarray}}
\def\BEQ{\begin{equation}}
\def\EEQ{\end{equation}}
\def\bref#1{(\ref{#1})}
\def\tmatrix#1#2#3#4{
    \left(\begin{array}{cc} #1 & #2 \\ #3 & #4 \end{array}\right)}
\def\LieSL{{\mbox{\bf SL}}}

\def\LieSU{{\mbox{\bf SU}}}
\def\LieB{{\mbox{\bf B}}}

\def\LieU{{\mbox{\bf U}}}
\def\Liesl{{\mbox{\bf sl}}}
\def\Liesu{{\mbox{\bf su}}}

\def\dprime{{\prime\prime}}

\def\vecr{{\mbox{\bf r}}}
\def\vecsigma{{\sigma\kern-2.25mm\sigma}}

\def\ta{\tilde{a}}

\def\tf{\tilde{f}}

\def\tw{\tilde{w}}
\def\tmu{\tilde{\mu}}
\def\tsigma{\tilde{\sigma}}

\def\tF{{\tilde{F}}}
\def\tH{{\tilde{H}}}

\def\tT{{\tilde{T}}}
\def\tU{{\tilde{U}}}

\def\FA{{\cal A}}
\def\FC{{\cal C}}

\def\FE{{\cal E}}
\def\FF{{\cal F}}
\def\FI{{\cal I}}

\def\FRL{{{\cal R}\kern1mm_\lambda}}
\def\FS{{\cal S}}

\def\hchi{\hat{\chi}}

\def\diff{{\rm d}}
\def\Ad{{\mbox{\rm Ad}}}
\def\res{{\mbox{\rm res}}}

\def\Sym{{\mbox{\rm Sym}}}
\def\Per{{\mbox{\rm Per}}}
\def\OAff{{\mbox{\rm OAff}}}

\def\twospace{{\Brr^{\!\lower2pt\hbox{\mbox{\rm\scriptsize{2}}}}}}
\def\threespace{{\Brr^{\!\lower2pt\hbox{\mbox{\rm\scriptsize{3}}}}}}
\def\cstar{{\Bcc^{\!\lower2pt\hbox{\mbox{\scriptsize{$\ast$}}}}}}
\def\Rplus{{\Brr^{\!\lower2pt\hbox{\mbox{\scriptsize{$+$}}}}}}
\def\Rplusnull{{\Brr_0^{\!\lower2pt\hbox{\mbox{\scriptsize{$+$}}}}}}
\def\Imag{{\mbox{\rm Im}}}
\def\id{{\mbox{\rm id}}}
\def\tr{{\mbox{\rm tr}}}

\def\ha{\hat{a}}
\def\hb{\hat{b}}
\def\hc{\hat{c}}
\def\hd{\hat{d}}

\def\hp{\hat{p}}
\def\hq{\hat{q}}

\def\hF{\hat{F}}
\def\hH{\hat{H}}

\def\hL{\hat{L}}

\def\hR{\hat{R}}
\def\hS{\hat{S}}
\def\hp{\hat{p}}

\def\hU{\hat{U}}
\def\hV{\hat{V}}

\def\hPsi{\hat{\Psi}}
\def\halpha{\hat{\alpha}}
\def\hbeta{\hat{\beta}}

\def\tT{\tilde{T}}
\def\alphanull{\alpha^\circ}
\def\betanull{\beta^\circ}
\def\pnull{p^\circ}
\def\qnull{q^\circ}

\def\Unull{{U^\circ}}
\def\Vnull{{V^\circ}}
\def\Fnull{F^\circ}
\def\Psinull{\Psi^\circ}

\def\zquer{{\overline{z}}}
\def\alphaquer{\overline{\alpha}}
\def\halphaquer{\overline{\halpha}}
\def\betaquer{\overline{\beta}}
\def\hbetaquer{\overline{\hbeta}}

\def\lambdaquer{\overline{\lambda}}
\def\nuquer{\overline{\nu}}
\def\rmin{{r_{\mbox{\scriptsize\rm min}}}}
\def\pder#1#2{{\partial #1\over\partial #2}}

\def\FCinvolution{\hat{\sigma}}
\def\FCstar{\FC^\ast}

\begin{document}

\renewcommand{\thefootnote}{\fnsymbol{footnote}}

\begin{center}
{\LARGE Discrete surfaces of constant mean curvature via dressing}

\vskip1cm
\begin{minipage}{6cm}
\begin{center}
G.~Haak\footnotemark[2]\\ Fachbereich Mathematik\\ TU Berlin\\ D-10623 Berlin
\end{center}
\end{minipage}
\vspace{0.5cm}

\end{center}

\footnotetext[1]{supported by Sonderforschungsbereich 288}

\section{Introduction}\refstepcounter{subsection} \label{introduction}
\message{[Introduction]}
The last ten years have seen large progress in the investigation of
surfaces of constant mean curvature. In particular, for CMC-surfaces
without umbilics many important new results were obtained. Starting
with the work of Wente~\cite{We:1}, CMC-tori were first constructed and
then classified in terms of algebraic geometric data
\cite{Ab:1,PiSt:1,Bo:1,ErKnTr:1,Ja:1}.

In the last three years, Bobenko and Pinkall introduced the notion of
discrete CMC-surfaces~\cite{BoPi:1,BoPi:2}. They showed, that these
surfaces, defined by elementary geometric properties, are related to
the solutions of an integrable discretization of the elliptic
$\sinh$-Gordon equation in the same way, as usual CMC-surfaces without
umbilics are related to solutions of the $\sinh$-Gordon equation itself.

This connection was established using the Lax representation of the
$\sinh$-Gordon equation, which is also at the starting point of the
dressing method for integrable systems.
Consequently, Pedit and Wu~\cite{PeWu:1} gave a method to construct discrete
CMC-surfaces, i.e., solutions of the discretized $\sinh$-Gordon
equation, by dressing of a discretized vacuum solution, a discrete cylinder.

For general conformal CMC-immersions, the investigation of symmetries
of the CMC-surface in $\threespace$ was started
in~\cite{DoHa:2}. There it was also shown, 
that, with the exception of cylinders,
there are no periodic surfaces in the dressing orbit of the standard cylinder.

On the other hand, in \cite{DoWu:1} it was shown, that there is a more
general definition of a dressing action, which gives all solutions of
finite type from the standard cylinder. In particular, all CMC-tori,
i.e., doubly periodic CMC-surfaces, are in this generalized dressing orbit.
In \cite{DoHa:3}, this was used to give an alternative derivation of
Pinkall and Sterling's classification of CMC tori.

Since the surfaces which were constructed by Pedit and Wu are all in
the dressing orbit of a discretized version of the standard cylinder,
it seemed natural to adapt the discussion of \cite{DoHa:2,DoHa:3} to discrete
CMC-surfaces. This is the goal of this paper.

We start in Section~\ref{disccyl} with the definition of the discrete
(standard) cylinder as an example of a discrete CMC-surface in the
sense of Bobenko and Pinkall~\cite{BoPi:1,BoPi:2} 
(Definition~\ref{discCMCdef}). This will also lead to the introduction
of so called extended frames for discrete CMC-surfaces.
These definitions will describe the
discrete analogues of CMC-surfaces without umbilics in an isothermic
parametrization.
Then we introduce the dressing action on the
cylinder to generate a large class of such discrete CMC-surfaces 
along the lines of \cite{PeWu:1}.

In Section~\ref{symofdisc}, we will introduce and investigate the
notion of a symmetry of a discrete CMC-surface. This definition is
modelled after the definition of translational symmetries for
continuous CMC-surfaces as given in \cite{DoHa:3}.
We will compute the transformation properties of several geometric
quantities related to the discrete CMC-surface under a symmetry.
The most important result here is Theorem~\ref{framesymmetry} which
gives the transformation equation of the extended frame of a
discrete CMC-surface.

In Section~\ref{symdressing} we will investigate the periodicity
conditions for discrete CMC-surfaces in more detail.
We will characterize all discrete,
periodic CMC-surfaces in the generalized dressing orbit of the
standard cylinder in terms of rational functions, 
similar to \cite[Theorems~3.6,3.7]{DoHa:3}.

In Section~\ref{algebrogeometric}, we will introduce a hyperelliptic
Riemann surface, which allows us, as in the continuous case, to
construct discrete periodic surfaces from algebro-geometric data.

\section{The $r$-dressing orbit of the discretized cylinder} 
\label{disccyl} \message{[disccyl]}

In this section we will define the discrete cylinder, the starting
point of our investigations, and its dressing orbit. This section
follows~\cite{PeWu:1} closely.

\newsection \label{cylinder}
Let us recall the formula for the standard cylinder
(see~\cite[Section~3.7]{DoHa:2}): The extended frame of the standard
cylinder is given by
\BEQ \label{standardcylinderframe}
F(z,\zquer,\lambda)=e^{(\lambda\inv
z-\lambda\zquer)A}=e^{x(\lambda\inv-\lambda)A+iy(\lambda\inv+\lambda)A},
\kern1cm\lambda\in S^1,
\EEQ
where we have used complex coordinates $z=x+iy$, $\zquer=x-iy$.
From this definition it follows, that
\BEA \label{contUnull}
U(x,y,\lambda) & = & F\inv\partial_xF=(\lambda\inv-\lambda)A \\
V(x,y,\lambda) & = & F\inv\partial_yF=i(\lambda+\lambda\inv)A.
\label{contVnull}
\EEA
The Euler method for the discretization of continuous dynamical
systems suggests the following discretization of the cylinder: Define
the matrices $\Unull,\Vnull\in\LieSU(2)_\sigma$ by
\BEA \label{discUnulldef}
\Unull(\lambda) & = & \frac{1}{\Delta_+}(I+r_1(\lambda\inv-\lambda)A), \\
\Vnull(\lambda) & = & \frac{1}{\Delta_-}(I+ir_2(\lambda\inv+\lambda)A),
\label{discVnulldef}
\EEA
where $r_1$ and $r_2$ are real, positive constants,
$r_1,r_2\in\Rplus$, and ($\lambda\in S^1$)
\BEA \label{Deltaplusdef}
\Delta_+ & = & \sqrt{\det(I+r_1(\lambda\inv-\lambda)A)}
=\sqrt{1-r_1^2(\lambda\inv-\lambda)^2}, \\
\Delta_- & = & \sqrt{\det(I+ir_2(\lambda\inv+\lambda)A)}
=\sqrt{1+r_2^2(\lambda\inv+\lambda)^2}. \label{Deltaminusdef}
\EEA
Here, we fix $\Delta_+$ and $\Delta_-$ by requiring, that they take
positive real values on $S^1$. This is possible, since for
$\lambda=e^{i\theta}\in S^1$ we have
\BEA \label{DeltaplusonS1}
\Delta_+ & = & \sqrt{1+4r_1^2\sin^2\theta}\in\Brr,\kern0.5cm\Delta_+\geq1,\\
\Delta_- & = & \sqrt{1+4r_2^2\cos^2\theta}\in\Brr,\kern0.5cm\Delta_-\geq1.
\label{DeltaminusonS1}
\EEA
From this we immediately get

\myprop{} {\em
For each pair $r_1,r_2\in\Rplus$, the functions $\Delta_+$ and
$\Delta_-$ defined above are even in $\lambda$ and real and positive
on $S^1$.
}

\separate
The functions $\Delta_+^2(\lambda)$ and $\Delta_-^2(\lambda)$ are
rational and can therefore be extended meromorphically to $\CPE$.

\mylemma{} {\em
For each real constant $r_1>0$, $\Delta_+^2$ has precisely four
simple zeores, which are located at $\lambda_+$, $-\lambda_+$,
$\lambda_+\inv$ and $-\lambda_+\inv$ on the real axis, where
\BEQ \label{lambdap}
\lambda_+
=\frac{1}{2r_1}+\sqrt{1+\frac{1}{4r_1^2}}>1.
\EEQ
For each real constant $r_2>0$, $\Delta_-^2$ has precisely four
simple zeores, which are located at $i\lambda_-$, $-i\lambda_-$,
$i\lambda_-\inv$ and $-i\lambda_-\inv$ on the imaginary axis, where
\BEQ \label{lambdam}
\lambda_-=\frac{1}{2r_2}+\sqrt{1+\frac{1}{4r_2^2}}>1,
\EEQ
In the limits we have
\BEA
\lim_{r_1\rightarrow0}\lambda_+=\infty,& & 
\lim_{r_1\rightarrow\infty}\lambda_+=1, \nonumber\\
\lim_{r_2\rightarrow0}\lambda_-=\infty,& & 
\lim_{r_2\rightarrow\infty}\lambda_-=1.\label{limits}
\EEA
}

\Proof Both $\Delta_+^2$ and $\Delta_-^2$ are rational with
numerators of degree $4$, i.e., they can have at most $4$ zeroes
with multiplicity.
If $\lambda_+\in\Bcc$ is a zero of
$\Delta_+^2=1-r_1^2(\lambda\inv-\lambda)^2$, then obviously also
$-\lambda_+$, $\lambda_+\inv$ and $-\lambda_+\inv$ are zeroes of
$\Delta_+^2$. If $\lambda_+$ is real and different from $1$, then
these zeroes are real and distinct. A direct calculation using
\BEQ \label{Deltapfactored}
\Delta_+^2=(1+r_1(\lambda\inv-\lambda))(1-r_1(\lambda\inv-\lambda))
\EEQ
shows, that
Eq.~\bref{lambdap} gives a zero of $\Delta_+^2$ on the real axis.
In the same way, the statement for $\Delta_-^2$ is proved using
\BEQ \label{Deltamfactored}
\Delta_-^2=(1+ir_2(\lambda\inv+\lambda))(1-ir_2(\lambda\inv+\lambda)).
\EEQ
The limits~\bref{limits} follow directly from Eqs.~\bref{lambdap}
and~\bref{lambdam}.
\QED

\newsection{} \label{loopgroups}
As in the continuous case (see e.g.~\cite[Section~2.2]{DoHa:3}),
we will interprete the $\lambda$-dependent matrices as taking values
in a certain loop group.

For each real constant $r$, $0<r<1$, 
let $\Lambda_r\LieSL(2,\Bcc)_\sigma$ denote the group of 
smooth maps $g(\lambda)$ from $C_r$, the circle of radius $r$, to
$\LieSL(2,\Bcc)$, which satisfy the twisting condition
\BEQ \label{DPWtwistcond}
g(-\lambda)=\sigma(g(\lambda)),
\EEQ
where $\sigma:\LieSL(2,\Bcc)\rightarrow\LieSL(2,\Bcc)$ 
is the group automorphism
of order $2$, which is given by conjugation with the Pauli matrix 
\BEQ
\sigma_3=\tmatrix100{-1}.
\EEQ
The Lie algebras of these groups, which we denote by
$\Lambda_r\Liesl(2,\Bcc)_\sigma$, consist of maps 
$x:C_r\rightarrow\Liesl(2,\Bcc)$, which satisfy a similar
twisting condition as the group elements
\BEQ
x(-\lambda)=\sigma_3 x(\lambda)\sigma_3.
\EEQ
In order to make these loop groups complex Banach Lie groups, we equip them,
as in \cite{DoPeWu:1}, with some $H^s$-topology for $s>{1\over2}$.
Elements of these twisted loop
groups are matrices with off-diagonal entries which are odd functions, and
diagonal entries which are even functions in the parameter $\lambda$.
All entries are in the Banach algebra $\FA_r$ of $H^s$-smooth
functions on $C_r$.

Furthermore, we will use the following subgroups of
$\Lambda_r\LieSL(2,\Bcc)_\sigma$: 
Let $\LieB$ be a subgroup of $\LieSL(2,\Bcc)$ and 
$\Lambda_{r,B}^+\LieSL(2,\Bcc)_\sigma$ be the group of maps in
$\Lambda_r\LieSL(2,\Bcc)_\sigma$, which can be extended to holomorphic maps on
\BEQ
I^{(r)}=\{\lambda\in\Bcc; |\lambda|<r\},
\EEQ
the interior of the circle $C_r$, and take values in $\LieB$ at $\lambda=0$.
Analogously, let $\Lambda_{r,B}^-\LieSL(2,\Bcc)_\sigma$ 
be the group of maps in $\Lambda_r\LieSL(2,\Bcc)_\sigma$, which can be extended
to the exterior 
\BEQ
E^{(r)}=\{\lambda\in\CPE;|\lambda|>r\}
\EEQ
of $C_r$ and take values in $\LieB$ at $\lambda=\infty$. 
If $\LieB=\{I\}$ (based loops) 
we write the subscript $\ast$ instead of $\LieB$, if
$\LieB=\LieSL(2,\Bcc)$ we omit the subscript for $\Lambda$ entirely.

Also, by an abuse of notation, we will denote by 
$\Lambda_r\LieSU(2)_\sigma$ the subgroup of maps in
$\Lambda_r\LieSL(2,\Bcc)_\sigma$, which can be extended holomorphically
to the open annulus
\BEQ
A^{(r)}=\{\lambda\in\Bcc;r<|\lambda|<\frac{1}{r}\}
\EEQ
and take values in $\LieSU(2)$ on the unit circle.
As usual, we will set
\BEQ
\Lambda\LieSU(2)_\sigma=\bigcup_{0<r<1}\Lambda_r\LieSU(2)_\sigma.
\EEQ
Corresponding to these subgroups, we analogously define Lie subalgebras of
$\Lambda_r\Liesl(2,\Bcc)_\sigma$.

We quote the following results from~\cite{Mc:1} and \cite{DoPeWu:1}:

\separate (i) For each solvable subgroup $\LieB$ of 
$\LieSL(2,\Bcc)$, which satisfies $\LieSU(2)\cdot\LieB=\LieSL(2,\Bcc)$ and
$\LieSU(2)\cap\LieB=\{I\}$, multiplication 
$$
\Lambda_r\LieSU(2)_\sigma\times\Lambda_{r,B}^+\LieSL(2,\Bcc)_\sigma
\longrightarrow\Lambda_r\LieSL(2,\Bcc)_\sigma
$$
is a diffeomorphism onto.
The associated splitting
\BEQ \label{Iwasawa}
g=F g_+
\EEQ
of an element $g$ of $\Lambda_r\LieSL(2,\Bcc)_\sigma$, s.t.\
$F\in\Lambda_r\LieSU(2)_\sigma$ and
$g_+\in\Lambda^+_{r,B}\LieSL(2,\Bcc)_\sigma$
will be called Iwasawa decomposition. In the following, we will fix the
group $\LieB$ as the group of upper triangular $2\times2$-matrices 
with real positive entries on the diagonal.

\separate (ii) Multiplication 
$$
\Lambda^-_{r,\ast}\LieSL(2,\Bcc)_\sigma\times\Lambda_r^+\LieSL(2,\Bcc)_\sigma
\longrightarrow\Lambda_r\LieSL(2,\Bcc)_\sigma
$$
is a diffeomorphism onto the open and dense subset 
$\Lambda^-_{r,\ast}\LieSL(2,\Bcc)_\sigma\cdot\Lambda_r^+\LieSL(2,\Bcc)_\sigma$
of $\Lambda_r\LieSL(2,\Bcc)_\sigma$, called the ``big cell'' \cite{SeWi:1}.
The associated splitting
\BEQ
g=g_-g_+
\EEQ
of an element $g$ of the big cell, where
$g_-\in\Lambda^-_{r,\ast}\LieSL(2,\Bcc)_\sigma$ and
$g_+\in\Lambda_r^+\LieSL(2,\Bcc)_\sigma$, 
will be called Birkhoff factorization.

\myprop{} {\em
Let $r_1,r_2\in\Rplus$ and define
\BEQ \label{rmindef}
0<\rmin(r_1,r_2)=\max(\lambda_+\inv,\lambda_-\inv)<1.
\EEQ
For each $\rmin(r_1,r_2)<r<1$, the matrices $\Unull$ and $\Vnull$
defined in~\bref{discUnulldef} and~\bref{discVnulldef} are in the twisted
loop group $\Lambda_r\LieSU(2)_\sigma$. In particular, they can be extended
holomorphically to
$A^{(\rmin)}=\{\lambda\in\CPE|\;\rmin<|\lambda|<\frac{1}{\rmin}\}$.
}

\Proof It is clear, that 
\BEQ \label{detUVnull}
\det \Unull(\lambda)=\det \Vnull(\lambda)=1.
\EEQ
By Proposition~\ref{cylinder}, $\Delta_+$ and 
$\Delta_-$ are real and do not vanish for $\lambda\in S^1$. Thus,
$\Unull$ and $\Vnull$ take values in $\LieSU(2)$ for $\lambda\in S^1$.
Since $\Delta_+\Unull$ and $\Delta_-\Vnull$ can be extended
holomorphically to $\cstar=\Bcc\setminus\{0\}$ and since, 
by Lemma~\ref{cylinder}, $\Delta_+^2$ and $\Delta_-^2$ have no zeroes
in $A^{(\rmin)}$, where $\rmin(r_1,r_2)$ is given by Eq.~\bref{rmindef}, 
we get, that $\Unull$ and $\Vnull$ can be extended
holomorphically to $A^{(\rmin)}$.
In addition, since $\Delta_+$ and $\Delta_-$ are even in $\lambda$, we
have, that $\Unull$ and $\Vnull$ satisfy the twisting
condition~\bref{DPWtwistcond}. Thus, for each $\rmin(r_1,r_2)<r<1$,
$\Unull$ and $\Vnull$ are in $\Lambda_r\LieSU(2)_\sigma$.
\QED

\newsection \label{Symsection}
We will define a map $\Fnull:\Bii^2\rightarrow\Lambda\LieSU(2)_\sigma$ by
\BEQ \label{Fnulldef}
\Fnull_{mn}(\lambda)=\Unull(\lambda)^m\Vnull(\lambda)^n.
\EEQ
By Proposition~\ref{loopgroups} and since $\Unull$ and $\Vnull$
commute, for all $\rmin(r_1,r_2)<r<1$, $\Fnull$ is a well defined map from
$\Bii^2$ to $\Lambda_r\LieSU(2)_\sigma$. Using Sym's formula, as in the
continuous case, we get a ``discrete surface''
$\Psinull_{mn}:\Bii^2\rightarrow\threespace$ in the spinor
representation $J:\threespace\rightarrow\Liesu(2)$,
$\vecr\mapsto-\frac{i}{2}\vecsigma\vecr$:
\BEQ \label{Sym}
J(\Psinull_{mn})=
\pder{}{\theta}|_{\theta=0}\Fnull_{mn}(\lambda=e^{i\theta})\cdot
\Fnull_{mn}(1)\inv+\frac{i}{2}\Fnull_{mn}(1)\sigma_3\Fnull_{mn}(1)\inv.
\EEQ

\newsection \label{discCMCdef}
For each $r_1,r_2\in\Rplus$, the 
discretized cylinder provides an example of a discrete
CMC-surface, as was shown in~\cite{PeWu:1}.
Let us define, what we understand by a discrete CMC-surface:

\mydefinition{} Let $\Psi_{mn}$ be a map from $\Bii^2$ into
$\threespace$. The map $\Psi_{mn}$ is
called a {\em discrete CMC-surface}\/ iff the following holds:
There exists a map $F_{mn}$ from $\Bii^2$ into
$\Lambda\LieSU(2)_\sigma$, s.t.\ for two real positive constants
$r_1$, $r_2$:
\BEA
U_{mn}(\lambda) & = & F_{mn}\inv F_{m+1,n}=
\frac{1}{\Delta_+}\tmatrix{\alpha_{mn}}{\lambda\inv r_1
p_{mn}-\lambda r_1 p_{mn}\inv}{\lambda\inv r_1
p_{mn}\inv-\lambda r_1 p_{mn}}{\alphaquer_{mn}}, \label{Uform}\\
V_{mn}(\lambda) & = & F_{mn}\inv F_{m,n+1}=
\frac{1}{\Delta_-}\tmatrix{\beta_{mn}}{i\lambda\inv r_2
q_{mn}+i\lambda r_2 q_{mn}\inv}{i\lambda\inv r_2
q_{mn}\inv+i\lambda r_2 q_{mn}}{\betaquer_{mn}},\label{Vform}
\EEA
where $p_{mn}$, $q_{mn}$ are positive real constants,
$\alpha_{mn}$ and $\beta_{mn}$
are complex constants, and $\Psi_{mn}$ is given in terms of $F_{mn}$
by Sym's formula~\bref{Sym}. The constants $r_1$, $r_2$ will be
called {\em lattice constants}.
The set of maps $F_{mn}:\Bii^2\rightarrow\Lambda\LieSU(2)_\sigma$, for which
$U_{mn}$ and $V_{mn}$ are of the form~\bref{Uform} and~\bref{Vform},
respectively, will be denoted
by $\FF(r_1,r_2)$. Its elements will be called {\em extended frames}.

\separate
Clearly, for each $r_1,r_2\in\Rplus$, the map $\Fnull$ is in
$\FF(r_1,r_2)$ with $p_{mn}=q_{mn}=1$. 
We call the associated CMC-immersion $\Psinull$ the
standard cylinder in $\FF(r_1,r_2)$.

\myremark{} 
It can easily be calculated (see~\bref{Ledge}), 
that the edge vectors $L^{\circ}_{mn}=J(\Psinull_{m+1,n}-\Psinull_{mn})$ are
independent of $m,n$. Thus, $\Psinull:\Bii^2\rightarrow\threespace$ is
a discrete analogue of a ruled surface which is generated by parallel
lines, i.e., a discrete analogue of a cylinder over a plane curve.
However, the reader should be aware that, unlike the continuous case,
for general lattice constants $r_1,r_2$, the standard cylinder will
not close up in one direction in the parameter space.

\separate
The motivation for Definition~\ref{discCMCdef} comes from the following

\mytheorem{} {\em
Let $\Psi_{mn}:\Bii^2\rightarrow\threespace$ be a discrete
CMC-surface with extended frame $F_{mn}$. Let us define $U_{mn}$ and
$V_{mn}$ as in~\bref{Uform} and~\bref{Vform}. Then for $\alpha_{mn}$ and
$\beta_{mn}$ we get
\BEA \label{alphabetrag}
|\alpha_{mn}|^2 & = & 1-r_1^2(p_{mn}-p_{mn}\inv)^2,\\
|\beta_{mn}|^2 & = & 1-r_2^2(q_{mn}-q_{mn}\inv)^2.
\label{betabetrag}
\EEA
Furthermore, $p=p_{mn}$, $q=q_{mn}$, $\alpha=\alpha_{mn}$,
$\beta=\beta_{mn}$, $p^\prime=p_{m,n+1}$,
$\alpha^\prime=\alpha_{m,n+1}$, $q^\prime=q_{m+1,n}$ and
$\beta^\prime=\beta_{m+1,n}$ satisfy
\BEA
pp^\prime=qq^\prime \label{closing}\\
\left(\frac{q}{p}+\frac{p}{q}\right)\alpha^\prime & = &
i\frac{r_1}{r_2}\left(\frac{p}{p^\prime}-\frac{p^\prime}{p}\right)
\betaquer+\left(\frac{q^\prime}{p}+\frac{p}{q^\prime}\right)\alphaquer,
\label{aprimeeq} \\
\left(\frac{q}{p}+\frac{p}{q}\right)\beta^\prime & = &
i\frac{r_2}{r_1}\left(\frac{q^\prime}{q}-\frac{q}{q^\prime}\right)
\alphaquer+\left(\frac{q^\prime}{p}+\frac{p}{q^\prime}\right)\betaquer,
\label{bprimeeq}
\EEA
and the discrete $\sinh$-Gordon equation
\BEQ \label{discSG}
\alpha^\prime\beta-\beta^\prime\alpha=ir_1r_2(p^\prime q+q^\prime p
-(p^\prime q)\inv-(q^\prime p)\inv).
\EEQ
}

\Proof
Eqs.~\bref{alphabetrag} and~\bref{betabetrag} follow from $\det
U_{mn}=\det V_{mn}=1$. Eqs.~\bref{closing}--\bref{discSG} are obtained
by a direct calculation from the integrability condition
\BEQ \label{discZCC}
U_{mn}V_{m+1,n}=V_{mn}U_{m,n+1}.
\EEQ
\QED

\myremark{}
1. By~\bref{closing} we can define
\BEA
p_{mn} & = & e^{-\frac12(\omega_{mn}+\omega_{m+1,n})}, \\
q_{mn} & = & e^{-\frac12(\omega_{mn}+\omega_{m,n+1})}.
\EEA
The equations \bref{aprimeeq},~\bref{bprimeeq}
and~\bref{discSG} then transform into those of Pedit and Wu
\cite[Theorem~4.1]{PeWu:1} for the real variables $\omega_{mn}$.
Bobenko and Pinkall first used the form of the Lax
operators~\bref{Uform},\bref{Vform} to
arrive at the integrable discretization~\bref{discSG} of the
$\sinh$-Gordon equation and to derive the elementary geometric
properties of the discretized CMC-surfaces.

2. It should be noted that in the discrete case there is no
natural definition of an associated family of discrete CMC-surfaces. 
The reason for
this loss of structure, compared to the continuous case, is the fact
that discrete CMC-surfaces are actually discrete isothermic
CMC-surfaces~\cite{BoPi:2}. 
In the continuous case the parametrization is isothermic only for real
Hopf differential (we only consider constant Hopf differential, i.e.,
no umbilics). 
If we choose $E\equiv 1$, as it was done in the definition of the
standard cylinder above, then these surfaces are the ones corresponding
to $\lambda=\pm1$ (see~the appendix of \cite{DoHa:1}).
Consequently, in the discrete case as we present it here, 
the Hopf differential is always normalized to $E\equiv1$, and the
associated family is reduced to the two surfaces for $\lambda=1$ and
$\lambda=-1$, which are congruent.

\mylemma{} {\em
For $r_1,r_2\in\Rplus$, let $F_{mn}\in\FF(r_1,r_2)$ be an extended frame.
Define $U_{mn}(\lambda)$
and $V_{mn}(\lambda)$ by~\bref{Uform},\bref{Vform} and
$\rmin(r_1,r_2)$ by~\bref{rmindef}. If
\BEQ \label{Finitial}
\mbox{\rm $F_{00}(\lambda)=I$ for all $\lambda\in S^1$,}
\EEQ
then for all $(m,n)\in\Bii^2$, 
\begin{enumerate}
\item[1.] $F_{mn}$, $U_{mn}$ and $V_{mn}$ can be extended holomorphically to
$A^{(\rmin)}$,
\item[2.] $\Delta_+U_{mn}$, $\Delta_-V_{mn}$ and
$\Delta_+^{|m|}\Delta_-^{|n|}F_{mn}$ can be extended holomorphically
to $\cstar$.
\end{enumerate}
}

\Proof
We prove first the second statement:
By~\bref{Uform} and~\bref{Vform}, $\Delta_+U_{mn}$ and
$\Delta_-V_{mn}$ can be continued holomorphically to $\cstar$ and the
same holds for $\Delta_+U_{mn}\inv$ and $\Delta_-V_{mn}\inv$. From
this and the initial condition~\bref{Finitial} it follows by induction
in $m$ and $n$, that also $\Delta_+^{|m|}\Delta_-^{|n|}F_{mn}$ can be
extended holomorphically to $\cstar$, which proves the second statement.

The first statement ist now a simple consequence of the fact, that
$\Delta_+$ and $\Delta_-$ have no zeroes on $A^{(\rmin)}$, i.e.,
$\Delta_+\inv$ and $\Delta_-\inv$ can be extended holomorphically to
$A^{(\rmin)}$.
\QED

\mydefinition{}
An extended frame $F_{mn}\in\FF(r_1,r_2)$ will be called {\em normalized},
if it satisfies~\bref{Finitial}. The subset of normalized extended frames
in $\FF(r_1,r_2)$ will be denoted by $\FF_0(r_1,r_2)$.

\separate
We will show in Section~\ref{framesymmetry}, 
that w.l.o.g.\ we can restrict our attention to
normalized extended frames.

Lemma~\ref{discCMCdef} shows, that for arbitrary 
$0<\rmin(r_1,r_2)<r<1$ and $(m,n)\in\Bii^2$,
$U_{mn}$ and $V_{mn}$ as well as the normalized extended frame
$F_{mn}(\lambda)$ are in $\Lambda_r\LieSU(2)_\sigma$.

\newsection{} \label{involutions}
For later use we introduce the following antiholomorphic involution
\BEQ
\tau:\lambda\mapsto\lambdaquer\inv.
\EEQ
Geometrically speaking, $\tau$ is the reflection at the unit circle in $\CPE$.
For a map $g(\lambda)$ from a subset of $\CPE$ to 
$\LieSL(2,\Bcc)$ we define
\BEQ
g^\ast(\lambda)=\overline{g(\tau(\lambda))}^\top.
\EEQ
Thus, if $F\in\Lambda_r\LieSL(2,\Bcc)_\sigma$ is defined and holomorphic on
$A^{(r)}$, then $F\in\Lambda_r\LieSU(2)_\sigma$ is equivalent to
\BEQ \label{realitycond}
F^\ast=F\inv.
\EEQ
For a scalar function $f(\lambda)$ we set
\BEQ
f^\ast(\lambda)=\overline{f(\tau(\lambda))}.
\EEQ
If $f$ is defined and holomorphic on a 
$\tau$-invariant neighbourhood of $S^1$, then
$f$ is real on $S^1$ iff $f^\ast=f$.

\newsection{} \label{dressing}
As in the continuous case, we can use the dressing action of the group
$\Lambda_r^+\LieSL(2,\Bcc)_\sigma$ on $\Lambda_r\LieSU(2)_\sigma$ to
generate new surfaces from old ones.

For $r_1,r_2\in\Rplus$, let 
$\Fnull:\Bii^2\rightarrow\Lambda_r\LieSU(2)_\sigma$ be an arbitrary
normalized extended frame in $\FF(r_1,r_2)$, not necessarily the discrete cylinder.
For arbitrary $\rmin(r_1,r_2)<r<1$ choose
$h_+\in\Lambda_r^+\LieSL(2,\Bcc)_\sigma$ and define a map
$F:\Bii^2\rightarrow\Lambda_r\LieSU(2)_\sigma$ by the Iwasawa splitting
\BEQ \label{Fdef}
h_+(\lambda)\Fnull_{mn}(\lambda)=F_{mn}(\lambda)p_+(m,n,\lambda),
\EEQ
where $p_+$ taking values in $\Lambda^+_r\LieSL(2,\Bcc)_\sigma$ is
chosen such that $F_{mn}(\lambda)=I$ for all $\lambda\in S^1$. We
have the following

\mytheorem{} {\em
Let $r_1,r_2\in\Rplus$ and $\Fnull\in\FF_0(r_1,r_2)$ be a normalized
extended frame.
Then for arbitrary $0<\rmin(r_1,r_2)<r<1$ and
$h_+\in\Lambda_r^+\LieSL(2,\Bcc)_\sigma$,
$F:\Bii^2\rightarrow\Lambda_r\LieSU(2)_\sigma$ defined
by~\bref{Fdef} is in $\FF_0(r_1,r_2)$.
}

\Proof
We define the following matrices
\BEA \label{discUdef}
U_{mn} & = & F_{mn}\inv F_{m+1,n}=p_+(m,n,\cdot)U^\circ_{mn}
p_+(m+1,n,\cdot)\inv,\\
V_{mn} & = & F_{mn}\inv F_{m,n+1}=p_+(m,n,\cdot)V^\circ_{mn}
p_+(m,n+1,\cdot)\inv,
\label{discVdef}
\EEA
where
\BEQ \label{Unullform}
U^\circ_{mn}=({\Fnull_{mn}})\inv\Fnull_{m+1,n}=
\frac{1}{\Delta_+}\tmatrix{\alphanull_{mn}}{\lambda\inv r_1
\pnull_{mn}-\lambda r_1(\pnull_{mn})\inv}{\lambda\inv r_1
(\pnull_{mn})\inv-\lambda r_1\pnull_{mn}}{\overline{\alphanull_{mn}}}
\EEQ
and
\BEQ \label{Vnullform}
V^\circ_{mn}=({\Fnull_{mn}})\inv\Fnull_{m,n+1}=
\frac{1}{\Delta_-}\tmatrix{\betanull_{mn}}{i\lambda\inv r_2
\qnull_{mn}+i\lambda r_2(\qnull_{mn})\inv}{i\lambda\inv r_2
(\qnull_{mn})\inv+i\lambda r_2\qnull_{mn}}{\overline{\betanull_{mn}}}.
\EEQ
Then, $U_{mn}$ and $V_{mn}$ are in $\Lambda_r\LieSU(2)_\sigma$.

Since $p_0(m,n)=p_+(m,n,\lambda=0)$ takes values
in the solvable Lie group $\LieB$, we can write
\BEQ \label{ppform}
p_0(m,n)=\tmatrix{e^{\omega_{mn}}}00{e^{-\omega_{mn}}},
\EEQ
where $\omega_{mn}\in\Brr$.
The matrix $\lambda\Delta_+U_{mn}$ is holomorphic in the vicinity of
$\lambda=0$ and takes the value 
\BEQ
p_0(m,n)\tmatrix0{r_1\pnull_{mn}}{r_1(\pnull_{mn})\inv}0p_0(m+1,n)\inv
\EEQ
at $\lambda=0$.
Thus, $U_{mn}$ is of the form
\BEQ \label{Uexpansion}
U_{mn}(\lambda)=\lambda\inv\frac{1}{\Delta_+}\tmatrix0{r_1
p_{mn}}{r_1 p_{mn}\inv}0+\tU_{mn}(\lambda),
\EEQ
with
\BEQ \label{pmndef}
p_{mn}=\pnull_{mn}e^{\omega_{mn}+\omega_{m+1,n}}\in\Rplus.
\EEQ
and $\tU_{mn}$ holomorphic in $\lambda$. Since
$U_{mn}\in\Lambda_r\LieSU(2)_\sigma$, we have
$\overline{U_{mn}(\lambda)}^\top=U_{mn}(\lambda)\inv$ for $\lambda\in
S^1$. Substituting this into~\bref{Uexpansion} and expanding
$\tU_{mn}$ as a power series in $\lambda$, we get
$\tU_{mn}(\lambda)=R_{mn}+S_{mn}\lambda$, where, by the twisting condition,
\BEQ
R_{mn}=\frac{1}{\Delta_+}\tmatrix{\alpha_{mn}}00{\alphaquer_{mn}}
\EEQ
for some complex constant $\alpha_{mn}$, and
\BEQ
S_{mn}=\frac{1}{\Delta_+}\tmatrix0{-r_1 p_{mn}\inv}{-r_1 p_{mn}}0.
\EEQ
This shows, that $U_{mn}$ is of the form~\bref{Uform}.
By a similar argument, Eq.~\bref{Vform}
follows with $\beta_{mn}\in\Bcc$ and
\BEQ \label{qmndef}
q_{mn}=\qnull_{mn}e^{\omega_{mn}+\omega_{m,n+1}}\in\Rplus.
\EEQ
Therefore, and since $F_{mn}$ satisfies~\bref{Finitial} by
construction, $F_{mn}$ is an extended frame.
\QED

\separate
Theorem~\ref{dressing} shows, that for arbitrary $r_1,r_2\in\Rplus$,
the groups $\Lambda_r^+\LieSL(2,\Bcc)_\sigma$, $\rmin(r_1,r_2)<r<1$
act on the set $\FF_0(r_1,r_2)$ of normalized extended frames with lattice
constants $r_1$ and $r_2$. We call this action the $r$-dressing action
on $\FF_0(r_1,r_2)$.

\section{Symmetries of discrete CMC-surfaces}
\label{symofdisc}\message{[symofdisc]}
\newsection \label{symmetrydef}
We will now define the symmetry group of a discrete CMC-surface:

\mydefinition{} Let $\Psi_{mn}:\Bii^2\rightarrow\threespace$ be a
discrete CMC-surface, then we denote by $\Sym(\Psi_{mn})$ the additive
group of all
pairs $(k,l)\in\Bii^2$, s.t.\
\BEQ \label{Symdef}
\Psi_{m+k,n+l}=\tT\Psi_{m,n}
\EEQ
for all $(m,n)\in\Bii^2$, where $\tT\in\OAff(\threespace)$, the set of
proper Euclidean motions of $\threespace$.
By $\Per(\Psi_{mn})$ we denote the subgroup of $\Sym(\Psi_{mn})$ of all pairs
$(k,l)\in\Bii^2$, for which Eq.~\bref{Symdef} is satisfied with
$\tT=\id$.
We will usually identify the elements of $\Sym(\Psi_{mn})$ and
$\Per(\Psi_{mn})$ with the corresponding translations in $\Bii^2$.

\separate
We will be mainly concerned with the transformation properties of
the extended frame $F_{mn}$ under translations in $\Sym(\Psi_{mn})$. A first
step in this direction is the following 

\mylemma{} {\em Let $\Psi_{mn}$ and $\hPsi_{mn}$ be two discrete 
CMC-surfaces with the same lattice constants $(r_1,r_2)$. 
Let $p_{mn}$, $q_{mn}$, $\hp_{mn}$ and $\hq_{mn}$ be the corresponding maps
from $\Bii^2$ to $\Rplus$, defined by~\bref{Uform},\bref{Vform}. If
$\Psi_{mn}$ and $\hPsi_{mn}$ differ only by a proper Euclidean motion, i.e., if
$\hPsi_{mn}=\tT\Psi_{mn}$ for all $(m,n)\in\Bii^2$, then
\BEQ
\hp_{mn}=p_{mn},\kern1cm\hq_{mn}=q_{mn},
\kern0.5cm\mbox{\rm for all $(m,n)\in\Bii^2$}.
\EEQ
}

\Proof Using Sym's formula and the Eqs.~\bref{Uform},\bref{Vform}
we get for the edge vectors
$L_{mn}=J(\Psi_{m+1,n}-\Psi_{mn})$ and
$R_{mn}=J(\Psi_{m,n+1}-\Psi_{mn})$ of $\Psi$
the following expressions:
\BEA
L_{mn} & = & 
-2ip_{mn}\Ad(F_{mn}(1))\tmatrix{r_1^2(p_{mn}-p_{mn}\inv)}{
r_1\alpha_{mn}}{r_1\alphaquer_{mn}}{-r_1^2(p_{mn}-p_{mn}\inv)}, \label{Ledge}\\
R_{mn} & = & 
-\frac{2i}{1+4r_2^2}q_{mn}\Ad(F_{mn}(1))\tmatrix{r_2^2(q_{mn}+q_{mn}\inv)}{
ir_2\beta_{mn}}{-ir_2\betaquer_{mn}}{-r_2^2(q_{mn}+q_{mn}\inv)}.
\label{Redge}
\EEA
From these equations we get with~\bref{alphabetrag}, \bref{betabetrag},
\BEA
\det(L_{mn})=4r_1^2 p_{mn}^2,\\
\det(R_{mn})=\frac{4r_2^2}{1+4r_2^2}q_{mn}^2.
\EEA
If $\hPsi_{mn}$ and $\Psi_{mn}$ differ only by a proper Euclidean motion, then
in the spinor representation we have
\BEQ
\hL_{mn}=U L_{mn}U\inv,\kern2cm\hR_{mn}=U R_{mn}U\inv,
\EEQ
for some unitary matrix $U$, where $\hL_{mn}$ and $\hR_{mn}$ are the edge
vectors of $\hPsi_{mn}$.
Therefore, $\det(\hL_{mn})=\det(L_{mn})$,
$\det(\hR_{mn})=\det(R_{mn})$, whence $\hp_{mn}^2=p_{mn}^2$,
$\hq_{mn}^2=q_{mn}^2$. Since $p_{mn}$, $q_{mn}$, $\hp_{mn}$ and
$\hq_{mn}$ take values in $\Rplus$, the claim follows.
\QED

\mycorollary{} {\em
Let $\Psi_{mn}$ be a discrete CMC-surface and let
$(k,l)\in\Sym(\Psi_{mn})$. If we define $p_{mn}$ and $q_{mn}$ by
Eqs.~\bref{Uform},\bref{Vform},
then $p_{m+k,n+l}=p_{mn}$ and $q_{m+k,n+l}=q_{mn}$.
}

\Proof If we set $\hPsi_{mn}=\Psi_{m+k,n+l}$, $\hp_{mn}=p_{m+k,n+l}$,
$\hq_{mn}=q_{m+k,n+l}$, then the proof follows
immediately from Lemma~\ref{symmetrydef}.
\QED

\newsection
From the proof of Lemma~\ref{symmetrydef} and the spinor
representation it is clear that
\BEQ
\tr(L_{mn}^2)=-2\det(L_{mn})=(4r_1 p_{mn})^2
\EEQ
and 
\BEQ
\tr(R_{mn}^2)=-2\det(R_{mn})=\left(\frac{4r_2}{\sqrt{1+4r_2^2}}q_{mn}\right)^2
\EEQ
are the squared lengths of the edge vectors of the discrete
CMC-surface. Therefore, for fixed $r_1$ and $r_2$, $p_{mn}$ and
$q_{mn}$ play the role of a discrete metric for the surface. We make the
following natural

\mydefinition{} Let $\Rplusnull$ be the set of non-negative real
numbers. The {\em metric of a discrete surface}
$\Psi_{mn}:\Bii^2\rightarrow\threespace$ is the map 
$g:\Bii^2\rightarrow\Rplusnull\times\Rplusnull$, defined by
\BEQ
g(m,n)=(|\hPsi_{m+1,n}-\hPsi_{mn}|,|\hPsi_{m,n+1}-\hPsi_{mn}|)
\EEQ
Two discrete surfaces
$\Psi_{mn},\hPsi_{mn}:\Bii^2\rightarrow\threespace$ will be called
{\em isometric} iff their metrics are the same, i.e., iff for all 
$(m,n)\in\Bii^2$,
\BEQ
|\hPsi_{m+1,n}-\hPsi_{mn}|=|\Psi_{m+1,n}-\Psi_{mn}|,\kern1cm
|\hPsi_{m,n+1}-\hPsi_{mn}|=|\Psi_{m,n+1}-\Psi_{mn}|.
\EEQ

\separate
Using this definition, Corollary~\ref{symmetrydef} above just states
(see~\cite[Corollary~2.6]{DoHa:2}) that elements of $\Sym(\Psi_{mn})$ act
as self-isometries of the surface.
We will therefore say, that the discrete surface $\Psi_{mn}$ has periodic
metric, if $\Sym(\Psi_{mn})$ contains a nontrivial element
$(k,l)\neq(0,0)$.

\newsection \label{framesymmetry}
Now we derive the transformation formula for an extended frame and the
Lax operators $U_{mn}$ and $V_{mn}$ under a symmetry transformation of
$\Psi_{mn}$. This will also show the ambiguity in the definition of an
extended frame for a given discrete CMC-surface.

\mylemma{} {\em
Let $r_1,r_2\in\Rplus$ and let $\Psi_{mn}$ and $\hPsi_{mn}$ be two
discrete CMC-surfaces with extended frames $F_{mn}$ and $\hF_{mn}$ in
$\FF(r_1,r_2)$, respectively. Let $U_{mn}$,
$V_{mn}$ be defined by Eqs.~\bref{Uform},\bref{Vform} and let
$\hU_{mn}$ and $\hV_{mn}$ be the analogous matrices for $\hF_{mn}$.
If $\Psi_{mn}$ and $\hPsi_{mn}$ are related by a proper Euclidean
motion $\tT$ on $\threespace$, i.e., if
$\hPsi_{mn}=\tT\Psi_{mn}$ for all $(m,n)\in\Bii^2$, then
\BEQ \label{chidef}
\hF_{mn}(\lambda)=\chi(\lambda)F_{mn}(\lambda)k^{(-1)^{m+n+1}},
\EEQ
where $\chi\in\Lambda\LieSU(2)_\sigma=\hF_{00}kF_{00}\inv$ and
$k\in\LieU(1)\subset\LieSU(2)$ is diagonal and $\lambda$-independent.

In addition, if $\tT(x)=R(x)+t$, where $R:\threespace\rightarrow\threespace$ 
is a rotation and $t\in\threespace$ a
translation, then
\BEQ \label{chitrafo}
J(R(x))=\chi(k,l,1)J(x)\chi(k,l,1)\inv,\kern1cm 
J(t)=\frac{i}{2}\pder{}{\theta}|_{\theta=0}\chi(k,l,\lambda=e^{i\theta})
\cdot\chi(k,l,1)\inv.
\EEQ
}

\Proof If $\Psi_{mn}$ and $\hPsi_{mn}$ are related by a proper
Euclidean motion, then we already know by Lemma~\ref{symmetrydef},
that $\hp_{mn}=p_{mn}$ and $\hq_{mn}=q_{mn}$ for all
$(m,n)\in\Bii^2$. This shows, with Eqs.~\bref{alphabetrag}
and~\bref{betabetrag}, that
\BEQ \label{abbetrag}
|\halpha_{mn}|=|\alpha_{mn}|,\kern2cm |\hbeta_{mn}|=|\beta_{mn}|.
\EEQ
If we define the edge vectors $L_{mn}$, $R_{mn}$, $\hL_{mn}$ and
$\hR_{mn}$ as in the proof of Lemma~\ref{symmetrydef}, then we know,
that the scalar products $\langle L_{mn},R_{mn}\rangle$ and
$\langle\hL_{mn},\hR_{mn}\rangle$ are equal for both surfaces. In the
spinor representation, this gives with Eqs.~\bref{Ledge} and~\bref{Redge}
\BEQ \label{LinR}
\halphaquer_{mn}\hbeta_{mn}-\halpha_{mn}\hbetaquer_{mn}
=\alphaquer_{mn}\beta_{mn}-\alpha_{mn}\betaquer_{mn}.
\EEQ
Together with~\bref{abbetrag} this yields
\BEQ
\halpha_{mn}=e^{i\phi_{mn}}\alpha_{mn},\kern2cm
\hbeta_{mn}=e^{i\phi_{mn}}\beta_{mn},
\EEQ
where $\phi_{mn}\in[0,2\pi)$. Now we invoke~\bref{aprimeeq}
and~\bref{bprimeeq}. This gives
\BEQ
\phi_{m+1,n}=\phi_{m,n+1}=-\phi_{mn},
\EEQ
which shows that
\BEQ
\halpha_{mn}=e^{i(-1)^{m+n}\phi}\alpha_{mn},\kern2cm
\hbeta_{mn}=e^{i(-1)^{m+n}\phi}\beta_{mn},
\EEQ
where we have set $\phi=\phi_{00}$.
By looking at~\bref{Uform},\bref{Vform}, we get
\BEQ
\hU_{mn}(\lambda)=k^{(-1)^{m+n}}U_{mn}k^{(-1)^{m+n}},
\EEQ
where 
\BEQ
k=\tmatrix{e^{i\frac{\phi}{2}}}00{e^{-i\frac{\phi}{2}}}.
\EEQ
Now, $U_{mn}$ and $V_{mn}$ ($\hU_{mn}$ and $\hV_{mn}$) determine
$F_{mn}$ ($\hF_{mn}$) up to left multiplication
with an element of $\Lambda\LieSU(2)_\sigma$.

Since $\tF_{mn}=F_{mn}k^{(-1)^{m+n+1}}$ and $\hF_{mn}$ both satisfy
\BEA
\tF_{mn}\inv\tF_{m+1,n} & = & \hF_{mn}\inv\hF_{m+1,n}=\hU_{mn},\nonumber\\
\tF_{mn}\inv\tF_{m,n+1} & = &
\hF_{mn}\inv\hF_{m,n+1}=\hV_{mn},
\EEA
we get Eq.~\bref{chidef}. Finally, Eq.~\bref{chitrafo} can be easily
derived from Sym's formula~\bref{Sym}.
\QED

\mycorollary{} {\em
Let $\Psi_{mn}:\Bii^2\rightarrow\threespace$ be a discrete CMC-surface
with lattice constants $r_1,r_2$. Let 
$\hF_{mn}(\lambda),F_{mn}(\lambda)\in\FF(r_1,r_2)$
be two extended frames for $\Psi_{mn}$. Then there exists
$k\in\LieU(1)$ and $\chi(\lambda)\in\Lambda\LieSU(2)_\sigma$, 
s.t.\
\BEQ
\hF_{mn}(\lambda)=\chi(\lambda)F_{mn}(\lambda)k^{(-1)^{m+n+1}},
\EEQ
where
\BEQ
\chi(1)=\pm I,\kern1cm\pder{}{\theta}|_{\theta=0}\chi(\lambda=e^{i\theta})=0.
\EEQ
}

\Proof If we set $\hPsi_{mn}=\Psi_{mn}$, 
then the proof follows immediately from Lemma~\ref{framesymmetry}.
\QED

\mytheorem{} {\em
Let $\Psi_{mn}:\Bii^2\rightarrow\threespace$ be a discrete CMC-surface with
extended frame $F_{mn}\in\FF(r_1,r_2)$, $r_1,r_2\in\Rplus$.
Then the following are equivalent:
\begin{description}
\item[1.] $(k,l)\in\Sym(\Psi_{mn})$,
\item[2.] the equation
\BEQ \label{Ftrafo}
F_{m+k,n+l}(\lambda)=\chi(k,l,\lambda)F_{mn}(\lambda)k^{(-1)^{m+n+1}},
\EEQ
holds, where $\chi(k,l,\lambda)\in\Lambda\LieSU(2)_\sigma$ can be
extended holomorphically to $A^{(\rmin)}$, and $k\in\LieU(1)$
is a diagonal $\lambda$-independent matrix.
\end{description}
If $(k,l)\in\Sym(\Psi_{mn})$, then $(k,l)\in\Per(\Psi_{mn})$ iff
\BEQ \label{oldpercond}
\chi(k,l,1)=\pm I\kern1cm\mbox{\rm and}\kern1cm
\pder{}{\theta}|_{\theta=0}\chi(k,l,\lambda=e^{i\theta})=0.
\EEQ
}

\Proof 
We set $\hPsi_{mn}=\Psi_{m+k,n+l}$ and
$\hF_{mn}(\lambda)=F_{m+k,n+l}(\lambda)$. Then 1.$\Rightarrow$2.\
follows from Lemma~\ref{framesymmetry}.

Conversely, if Eq.~\bref{Ftrafo} is satisfied with $\chi$ and
$k_{mn}$ as above, then by Sym's formula~\bref{Sym} we have
\BEQ \label{Symtrafo}
\Psi_{m+k,n+l}=\chi(k,l,1)\Psi_{mn}\chi(k,l,1)\inv
+\frac{i}{2}\pder{}{\theta}|_{\theta=0}\chi(k,l,\lambda=e^{i\theta})
\cdot\chi(k,l,1)\inv.
\EEQ
Thus, $\Psi_{m+k,n+l}$ and $\Psi_{mn}$ differ by a proper Euclidean
motion $\tT$, which can be written as
\BEQ
\tT(x)=R(x)+t,\kern1cm\mbox{\rm for all $x\in\threespace$},
\EEQ
where $R:\threespace\rightarrow\threespace$ 
is a rotation and $t\in\threespace$ describes a
translation. In the spinor representation, $R$ and $t$ are given
by
\BEQ
J(R(x))=\chi(k,l,1)J(x)\chi(k,l,1)\inv,\kern1cm 
J(t)=\frac{i}{2}\pder{}{\theta}|_{\theta=0}\chi(k,l,\lambda=e^{i\theta})
\cdot\chi(k,l,1)\inv.
\EEQ
From this, 1.\ and the rest of the statement follow.
\QED

\separate
From now on we normalize, analogous to the continuous case,
the extended frames for the discrete CMC-surfaces by~\bref{Finitial}.
Up to proper Euclidean motions we still get,
by Lemma~\ref{framesymmetry}, all discrete CMC-surfaces in this way.

\section{Symmetric surfaces in the $r$-dressing orbit of the
cylinder}\label{symdressing}\message{[symdressing]}

In this section for fixed lattice constants $r_1,r_2\in\Rplus$ and
fixed $(k,l)\in\Bii^2\setminus\{(0,0)\}$,
we want to classify those discrete CMC-surfaces
$\Psi_{mn}:\Bii^2\rightarrow\threespace$ in the
dressing orbit of the cylinder, for which $(k,l)\in\Sym(\Psi_{mn})$.

Here, if we talk about the dressing orbit without
specifying a radius $r$, we always mean the orbit under the action of
the union
$$
\bigcup_{\rmin(r_1,r_2)<r<1}\Lambda_r^+\LieSL(2,\Bcc)_\sigma,
$$
where $\rmin(r_1,r_2)$ was defined in~\bref{rmindef}.

\newsection
Let $F_{mn}(\lambda)\in\FF_0(r_1,r_2)$ be defined by the dressing action of
$h_+\in\Lambda^+_r\LieSL(2,\Bcc)_\sigma$, $\rmin(r_1,r_2)<r<1$,
on the extended frame of the cylinder with lattice
constants $r_1,r_2\in\Rplus$. I.e.,
\BEQ \label{Fdressingdef}
h_+(\lambda)\Unull(\lambda)^m\Vnull(\lambda)^n
=F_{mn}(\lambda)p_+(n,m,\lambda),\kern1cm 
p_+(z,\lambda)\in\Lambda^+_{r,B}\LieSL(2,\Bcc)_\sigma.
\EEQ
Under the translation $(m,n)\mapsto(m+k,n+l)$, $(k,l)\in\Bii^2$, the
extended frame $F_{mn}(\lambda)$ transforms like
\BEQ
F_{m+k,n+l}(\lambda)=Q(\lambda)F_{mn}(\lambda)r_+(m,n,\lambda),
\EEQ
where
\BEQ
Q(\lambda)=h_+(\Unull)^k(\Vnull)^lh_+\inv
\EEQ
and
\BEQ
r_+(m,n,\lambda)=p_+(m,n,\lambda)p_+(m+k,n+l,\lambda)\inv.
\EEQ

\newsection \label{chiform}
By Theorem~\ref{framesymmetry} we have $(k,l)\in\Sym(\Psi_{mn})$
iff
\BEQ \label{Qchirel}
Q(\lambda)F_{mn}(\lambda)r_+(m,n,\lambda)=\chi(\lambda)F_{mn}(\lambda)k_{mn},
\EEQ
where $k_{mn}\in\LieSU(2)$ is a diagonal matrix,
$\chi(\lambda)\in\Lambda\LieSU(2,\Bcc)_\sigma$, and by Lemma~\ref{discCMCdef},
\BEQ
\Delta_+^{|k|}\Delta_-^{|l|}\chi(\lambda)=
\Delta_+^{|k|}\Delta_-^{|l|}F_{kl}(\lambda)k_{00}\inv
\EEQ
can be extended holomorphically to $\cstar$. In the following we
will derive further conditions on the matrix $\chi(\lambda)$.

The initial condition~\bref{Finitial} together with Eq.~\bref{Qchirel}
implies
\BEQ \label{chiQrel}
\chi=Q\cdot R_+,
\EEQ
with
\BEQ
R_+(\lambda)=r_+(0,0,\lambda)k_{00}\inv
\in\Lambda_r^+\LieSL(2,\Bcc)_\sigma.
\EEQ
Thus, $F_{mn}(\lambda)$ is invariant under the $r$-dressing
transformation with $R_+$,
\BEQ \label{thisequation}
R_+(\lambda)F_{mn}(\lambda)=F_{mn}(\lambda)r_+(m,n,\lambda)k_{mn}\inv.
\EEQ

\mylemma{} {\em
The matrix $\chi(\lambda)$ is of the form
\BEQ \label{chiHplus}
\chi=h_+Hh_+\inv,
\EEQ
where
\BEQ \label{Hdef}
H=(\Unull)^k(\Vnull)^l w_+
\EEQ
and $w_+=h_+\inv R_+h_+$ commute with $A=\tmatrix0110$.
}

\Proof
Substituting~\bref{Fdressingdef} into~\bref{thisequation} and rearranging
terms gives
\BEQ
\Vnull^{-n}\Unull^{-m}h_+\inv
R_+h_+\Unull^m\Vnull^n=
p_+(m,n,\lambda)\inv r_+(m,n,\lambda)k_{mn}\inv 
p_+(m,n,\lambda).
\EEQ
Abbreviating
\BEQ
V_+(m,n,\lambda)=p_+(m,n,\lambda)\inv
r_+(m,n,\lambda)k_{mn}\inv p_+(m,n,\lambda),
\EEQ
and using the definition of $w_+$, this is
\BEQ \label{wpVp}
\Vnull^{-n}\Unull^{-m}w_+\Unull^m
\Vnull^n=V_+(m,n,\lambda).
\EEQ
This yields that
\BEQ \label{Zinplus}
Z(m,n,\lambda)=\Vnull^{-n}\Unull^{-m}w_+\Unull^m\Vnull^n
\in\Lambda_r^+\LieSL(2,\Bcc)_\sigma,
\EEQ
in particular, that $Z(m,n,\lambda)$ is holomorphic on $I^{(r)}$, $\rmin<r<1$,
for all $m,n\in\Bii$.
Let $D=\frac{1}{\sqrt2}\tmatrix11{-1}1$, then
\BEQ
DAD\inv=\tmatrix100{-1}.
\EEQ
If we set
\BEQ
\tw_+=D w_+ D\inv=\tmatrix{\tw_a}{\tw_b}{\tw_c}{\tw_d}
\EEQ
then
\BEA
\lefteqn{DZ(m,n,\lambda)D\inv=} & & \nonumber\\
& =\tmatrix{\tw_a}{
\tw_b\left(\frac{1+r_1(\lambda\inv-\lambda)}{
1-r_1(\lambda\inv-\lambda)}\right)^m
\left(\frac{1+ir_2(\lambda\inv+\lambda)}{
1-ir_2(\lambda\inv+\lambda)}\right)^n}{
\tw_c\left(\frac{1-r_1(\lambda\inv-\lambda)}{
1+r_1(\lambda\inv-\lambda)}\right)^m
\left(\frac{1-ir_2(\lambda\inv+\lambda)}{
1+ir_2(\lambda\inv+\lambda)}\right)^n}{\tw_d}.& \label{tweq}
\EEA
If we define
\BEQ
S_{r_1}(\lambda)=\frac{1+r_1(\lambda\inv-\lambda)}{
1-r_1(\lambda\inv-\lambda)}
\EEQ
and
\BEQ
T_{r_2}(\lambda)=\frac{1+ir_2(\lambda\inv+\lambda)}{
1-ir_2(\lambda\inv+\lambda)},
\EEQ
then the off-diagonal entries of $Z(m,n,\lambda)$ are of the form
\BEQ
\tw_bS_{r_1}(\lambda)^mT_{r_2}(\lambda)^n
\EEQ
and
\BEQ
\tw_cS_{r_1}(\lambda)^{-m}T_{r_2}(\lambda)^{-n}.
\EEQ
Using~\bref{Deltapfactored} and the notation of
Lemma~\ref{cylinder}, we see that $S_{r_1}$ has two simple
zeroes at $\lambda_+$ and $-\lambda_+\inv$, and two simple poles at
$-\lambda_+$ and $\lambda_+\inv$. Analogously,
using~\bref{Deltamfactored}, we get that $T_{r_2}$ has two simple zeroes
at $i\lambda_-$ and $-i\lambda_-\inv$ and two simple poles at
$-i\lambda_-$ and $i\lambda_-\inv$.

Assume now, that $\tw_b$ and $\tw_c$ do not vanish identically. 
Then, for $m\neq0$ and $n\neq0$, the r.h.s.\ of~\bref{tweq}
has a simple pole on the circle $C_\rmin\subset I^{(r)}$. This
contradicts
$DZ(m,n,\lambda)D\inv\in\Lambda_r^+\LieSL(2,\Bcc)_\sigma$. Therefore,
$\tw_b=\tw_c\equiv0$ and $\tw_+$ is diagonal and commutes with
$\sigma_3=DAD\inv$. Thus, $w_+$ commutes with $A$.

This implies
\BEQ
V_+(m,n,\lambda)=V_+(0,0,\lambda)=w_+(\lambda).
\EEQ
The definition of $w_+$ gives
\BEQ \label{Rpwp}
R_+=h_+w_+h_+\inv.
\EEQ
Substituting~\bref{Rpwp} into~\bref{chiQrel} finally
gives~\bref{chiHplus}, with $H$ defined by~\bref{Hdef}. Since $\Unull$,
$\Vnull$ and $w_+$ commute with $A$, also $H$ does.
\QED

\newsection{} \label{Ddigression}
In the proof of Lemma~\ref{chiform} we have
diagonalized $A$. We will use this method to derive a more convenient
expression for the matrix $H$ introduced in the lemma.

Recall, that $\FA_r$ was the Banach algebra defining the topology of
the loop groups introduced in Section~\ref{loopgroups}.
Let $\FA_r^+$ be the subalgebra of $\FA_r$ which consists of those
functions which can be continued holomorphically to $I^{(r)}$. 
Lemma~\ref{chiform} together with $\det w_+=1$ implies
(see~\cite[Section~3.3]{DoHa:3}), that the matrix
$w_+$ is of the form
\BEQ
w_+=e^{f_+A}=\tmatrix{\cosh(f_+)}{\sinh(f_+)}{\sinh(f_+)}{\cosh(f_+)},
\EEQ
where $f_+\in\FA_r^+$ is odd in $\lambda$. 
For $\tw_+=Dw_+D\inv$ this amounts to
$\tw_+=e^{f_+\sigma_3}$. 
The matrix $\tH=DHD\inv=D(\Unull)^k(\Vnull)^lw_+D\inv$ is of the form
\BEQ
\tH=\frac{1}{\Delta_+^k\Delta_-^l}\tmatrix{r_+}00{r_-},
\EEQ
where
\BEA
r_+(\lambda) & = & (1+r_1(\lambda\inv-\lambda))^k
(1+ir_2(\lambda\inv+\lambda))^l e^{f_+},\\
r_-(\lambda) & = & (1-r_1(\lambda\inv-\lambda))^k
(1-ir_2(\lambda\inv+\lambda))^l e^{-f_+}.
\EEA
A short computation using $\det\tH=1$ and the parity of $f_+$
shows that we can rewrite this as
\BEQ \label{tHform}
\tH=\frac{1}{\Delta_+^{|k|}\Delta_-^{|l|}}
\tmatrix{\hp(\lambda)}00{\hp(-\lambda)},
\EEQ
where
\BEQ \label{hpdef}
\hp(\lambda)=(1+\epsilon_kr_1(\lambda\inv-\lambda))^{|k|}
(1+i\epsilon_lr_2(\lambda\inv+\lambda))^{|l|}e^{f_+},
\EEQ
with $k=\epsilon_k|k|$, $l=\epsilon_l|l|$. Obviously,
\BEQ \label{hpinpunctured}
\mbox{\rm $\hp(\lambda)$ is defined and holomorphic in $I^{(r)}\setminus\{0\}$.}
\EEQ
For $H$ this gives

\mylemma{} {\em
The matrix $H$ introduced in Lemma~\ref{Qchirel} is of the form
\BEQ \label{alphabetadef}
H=\alpha I+\beta A,\kern1cm\alpha,\beta\in\FA_r,
\EEQ
where
\BEQ \label{a2b2}
\alpha^2-\beta^2=1,
\EEQ
\BEQ \label{alphabetaparity}
\mbox{\rm $\alpha$ is even in $\lambda$ and $\beta$ is odd in
$\lambda$}.
\EEQ
and
\BEQ \label{deusorigin}
\alpha+\beta=\frac{1}{\Delta_+^{|k|}\Delta_-^{|l|}}\hp,
\EEQ
with $\hp$ given by~\bref{hpdef}.
}

\Proof
Eq.~\bref{alphabetadef} follows from~$[H,A]=0$, since $A$ is regular
semisimple. Then~$\det H=\alpha^2-\beta^2=1$ gives~\bref{a2b2},
and~\bref{alphabetaparity} follows from the twisting
condition~\bref{DPWtwistcond}. Finally,~\bref{deusorigin} follows 
from~\bref{tHform} by conjugation with $D\inv$.
\QED

\newsection{} \label{alphabetaintro}
Let us derive further properties of the functions $\alpha$ and $\beta$
introduced in Lemma~\ref{Ddigression}.
Eq.~\bref{alphabetadef} together with~\bref{a2b2} yields
\BEQ \label{Hinv}
H\inv=\alpha I-\beta A.
\EEQ
Since $(k,l)\neq(0,0)$, we get from~\bref{Hdef} and Lemma~\ref{cylinder}, 
that $H$ has a pole in $I^{(r)}$, $r>\rmin(r_1,r_2)$. Therefore,
by~\bref{a2b2},
\BEQ \label{betanontrivial}
\beta\not\equiv0.
\EEQ
Let us also define the matrix
\BEQ \label{hHdef}
\hH=\Delta_+^{|k|}\Delta_-^{|l|}H
=\Delta_+^{|k|}\Delta_-^{|l|}(\Unull)^k(\Vnull)^lw_+=\halpha I+\hbeta A,
\EEQ
with
\BEQ \label{hahbdef}
\halpha=\Delta_+^{|k|}\Delta_-^{|l|}
\alpha,\kern1cm
\hbeta=\Delta_+^{|k|}\Delta_-^{|l|}\beta.
\EEQ
and
\BEQ \label{halphahbetarel}
\det\hH=\halpha^2-\hbeta^2=\Delta_+^{2|k|}\Delta_-^{2|l|}=
(1-r_1^2(\lambda\inv-\lambda)^2)^{|k|}(1+r_2^2(\lambda\inv+\lambda)^2)^{|l|}.
\EEQ
Since $\Delta_+$ and $\Delta_-$ are even in $\lambda$ we have
\BEQ \label{halphahbetaparity}
\mbox{\rm$\halpha$ is even in $\lambda$ and $\hbeta$ is odd in
$\lambda$}.
\EEQ
Since by Theorem~\ref{framesymmetry},
\BEQ \label{hchihH}
\hchi(\lambda)=\Delta_+^{|k|}\Delta_-^{|l|}\chi(\lambda)=h_+\hH h_+\inv
\EEQ
is holomorphic on $\cstar$, taking the trace of $\hchi$ and $\hchi^2$
shows that $\halpha$ and $\hbeta^2$ have holomorphic extensions to
$\cstar$.
Now, since $w_+\in\Lambda^+_r\LieSL(2,\Bcc)_\sigma$ and since
$\Delta_+\Unull$ and $\Delta_-\Vnull$ as well as
$\Delta_+(\Unull)\inv$ and $\Delta_-(\Vnull)\inv$ are holomorphic on
$\cstar$, with meromorphic extension to $\lambda=0$ and $\lambda=\infty$,
also $\hH(\lambda)$ is meromorphic on
$I^{(r)}$. Thus $\hchi$ has a meromorphic extension to $\lambda=0$,
and by unitarity also to $\lambda=\infty$. This yields
\BEQ \label{halphahbetarational}
\mbox{\rm $\halpha$ and $\hbeta^2$ are rational functions which are
holomorphic on $\cstar$.}
\EEQ
From~\bref{deusorigin}, we get
\BEQ \label{deus}
\halpha+\hbeta=\hp,
\EEQ
where $\hp(\lambda)$ is defined by~\bref{hpdef}. 
From this and~\bref{halphahbetarational}
it follows, that
\BEQ
\mbox{\rm $\hbeta$ has a meromorphic extension to $I^{(r)}$.}
\EEQ

\newsection \label{Sdef}
We now define the matrices
\BEQ \label{Sdefeq}
S=h_+Ah_+\inv=\tmatrix abcd,
\EEQ
and
\BEQ
\hS=\hbeta S=\tmatrix{\ha}{\hb}{\hc}{\hd}.
\EEQ
With~\bref{alphabetadef} and~\bref{hchihH}, we have
\BEQ
\chi=\alpha I+\beta S
\EEQ
and
\BEQ \label{hchihahS}
\hchi=\halpha I+\hbeta S=\halpha I+\hS.
\EEQ
Using~\bref{Hinv} we get
\BEQ \label{chiinv}
\chi\inv=\alpha I-\beta S.
\EEQ
Since $\hchi$ and $\halpha$ are meromorphic on $\CPE$ and holomorphic
on $\cstar$, we get from~\bref{hchihahS}
\BEQ \label{hSentriesrational}
\mbox{\rm $\ha,\hb,\hc,\hd$ are rational and holomorphic on $\cstar$.}
\EEQ
Clearly, we have $\tr S=\tr\hS=0$, whence
\BEQ \label{deq}
d=-a,\kern3cm\hd=-\ha.
\EEQ
Also, in view of~\bref{halphahbetaparity}, the twisting condition for
$\Lambda_r\LieSL(2,\Bcc)_\sigma$ implies
\BEQ \label{abparity}
\mbox{\rm $\ha,b,c$ are even in $\lambda$, $a,\hb,\hc$ are odd in
$\lambda$.}
\EEQ
Since $S^2=I$, we get with~\bref{deq}
\BEQ \label{abcrel}
a^2+bc=1
\EEQ
Since $\hS=\hbeta S$, we also have
\BEQ \label{hahbhbeta}
\ha^2+\hb\hc=\hbeta^2.
\EEQ
The unitarity of $\chi(\lambda)$ on $S^1$ is in view of
Section~\ref{involutions} and~\bref{chiinv} equivalent with
\BEQ \label{alphareal}
\alpha^\ast=\alpha,
\EEQ
\BEQ \label{hSreality}
\hS^\ast=-\hS.
\EEQ
In particular,~\bref{hSreality} is equivalent with
\BEQ \label{hahbreal}
\ha^\ast=-\ha,\kern3cm\hb^\ast=-\hc.
\EEQ
By Proposition~\ref{cylinder}, we see
that~\bref{alphareal} is equivalent to
\BEQ \label{halphareal}
\halpha^\ast=\halpha,
\EEQ
Now we consider the squares of $\halpha,\hbeta,a,b,c$ and
$\ha,\hb,\hc$.
First we note
\BEQ \label{alpha2beta2real}
\mbox{\rm $\halpha^2$ and $\hbeta^2$ are real on $S^1$.}
\EEQ
This follows for $\halpha^2$ from~\bref{halphareal} and for $\hbeta^2$
from~\bref{halphahbetarel} and Proposition~\ref{cylinder}.
Next,~\bref{hahbreal} implies
\BEQ \label{ha2nonpositive}
\mbox{\rm $\ha^2$ is non-positive on $S^1$.}
\EEQ
Substituting this and~\bref{hahbreal} into~\bref{hahbhbeta} gives
\BEQ \label{b2nonpositive}
\mbox{\rm $\hbeta^2$ is non-positive on $S^1$.}
\EEQ
Since $\ha=\hbeta a$, we know $\ha^2=\hbeta^2a^2$. In particular,
by~\bref{hSentriesrational}, \bref{halphahbetarational},
\bref{ha2nonpositive} and~\bref{b2nonpositive},
\BEQ \label{a2rational}
\mbox{\rm $a^2$ is a rational function, real and non-negative on
$S^1$.}
\EEQ
Since $a^2$ is by definition also holomorphic at $\lambda=0$,
it follows that 
\BEQ \label{a2at0}
\mbox{\rm $a^2$ is locally holomorphic around $0$ and $\infty$.}
\EEQ
For $b^2$ and $c^2$ one argues
similarly. E.g.~$b^2=\frac{\hb^2}{\hbeta^2}$ is clearly meromorphic on
$\cstar$ and is also, by the definition of $b$, holomorphic at
$\lambda=0$. From~\bref{hahbreal} we obtain that
$(b^2)^\ast=\frac{(\hb^2)^\ast}{\hbeta^2}=\frac{\hc^2}{\hbeta^2}=c^2$ is
also holomorphic at $\lambda=0$. This shows, that $b^2$ is meromorphic
on $\CPE$ and thus rational. Altogether we have shown
\BEQ \label{b2c2rational}
\mbox{\rm $b^2$ and $c^2$ are rational and finite at $0$ and
$\infty$,}
\EEQ
and
\BEQ \label{b2c2rel}
(b^2)^\ast=c^2.
\EEQ
Next, from~\bref{Sdefeq} we see that $a(\lambda=0)=0$ and
$b(\lambda=0)=c(\lambda=0)\inv$. Since $a$ is odd in $\lambda$, we obtain
\BEQ \label{azero}
\mbox{\rm $a^2$ has a zero of order $2(2n-1)$ for some $n>0$ at
$\lambda=0$,}
\EEQ
\BEQ \label{bc1}
b(\lambda=0)c(\lambda=0)=1.
\EEQ
We also note that the relations 
\BEQ
\ha^2=\hbeta^2a^2,\kern1cm\hb^2=\hbeta^2b^2,\kern1cm\hc^2=\hbeta^2c^2
\EEQ
show that $a^2$, $b^2$, and $c^2$ can have poles only where 
$\hbeta^2$ has a zero.

Finally, from~\bref{hSreality} we obtain $(\hbeta b)^\ast=-(\hbeta
c)$. Hence~\bref{abcrel} implies
\BEQ
\hbeta^2=\hbeta^2a^2+\hbeta b\cdot\hbeta c=\hbeta^2a^2-(\hbeta b)(\hbeta
b)^\ast.
\EEQ
Therefore, on $S^1$ we obtain $\hbeta^2(a^2-1)=|\hbeta b|^2$. Since
$\hbeta^2$ is non-positive on $S^1$ by~\bref{b2nonpositive}, and
$\hbeta^2\not\equiv0$ by~\bref{betanontrivial}, we have
$a^2-1\leq0$. Thus,
\BEQ \label{aonS}
0\leq a^2(\lambda)\leq 1 \kern1cm\mbox{\rm for $\lambda\in S^1$}.
\EEQ

\newsection{} \label{lastnecc}
In the last section we considered the matrix $S=\tmatrix abcd$ and we
listed properties of $a,b,c,d$. In the rest of this paper we will
characterize $\Sym(\Psi_{mn})$ in terms of $a,b,c,d$.
The entries of elements of $\Lambda_r^+\LieSL(2,\Bcc)_\sigma$ are elements
of $\FA_r^+$. To be sure, that
to such $a,b,c,d$ there exists an $h_+$ satisfying~\bref{Sdefeq},
producing $F$ ---and thus $\Psi$--- for which
$(k,l)\in\Sym(\Psi_{mn})$, we slightly extend Theorem~3.6 from \cite{DoHa:3}:

\mytheorem{} {\em Let $a,b,c,d\in\FA_R^+$, $0<R<1$, where $a,d$ are odd and
$b,c$ are even. Then $S=\tmatrix abcd$ is of the form $S=h_+Ah_+\inv$
for some $0<r\leq R$ and $h_+\in\Lambda_r^+\LieSL(2,\Bcc)_\sigma$ iff
\BEQ
d=-a,
\EEQ
\BEQ \label{abcinth}
a^2+bc=1,
\EEQ
\BEQ
b(\lambda=0)\neq0.
\EEQ
Furthermore, if for some $0<r^\prime<R<1$, $b$ is the square of a
holomorphic function on the closure $I^{(r^\prime)}$, then we can find
such an $h_+$ for some $r\geq r^\prime$.
}

\Proof
All but the last statement are taken from Theorem~3.6 in~\cite{DoHa:3}.
The last statement follows from the special choice
\BEQ
h_+=\frac{1}{\sqrt{c}}\tmatrix1a0c
\EEQ
given in the proof of this Theorem.
\QED

\newsection{} \label{nectheorem}
Let us collect the necessary conditions we have derived in
Sections~\ref{chiform}--\ref{Sdef}:

\mytheorem{} {\em
Let $\Psi_{mn}:\Bii^2\rightarrow\threespace$ be a discrete CMC-surface
with extended frame $F_{mn}(\lambda)\in\FF_0(r_1,r_2)$, s.t.\
$F_{mn}(\lambda)$ is given by dressing the 
cylinder under the $r$-dressing~\bref{Fdressingdef} with some 
$h_+\in\Lambda_r^+\LieSL(2,\Bcc)_\sigma$, $\rmin(r_1,r_2)<r<1$.
Assume also, that for
$(k,l)\in\Bii^2$, $(k,l)\neq(0,0)$,
$F_{m+k,n+l}(\lambda)=\chi(\lambda)F_{mn}(\lambda)$,
i.e., $(k,l)\in\Sym(\Psi_{mn})$.
Define $h_+Ah_+\inv=\tmatrix abcd$. Then $d=-a$
and the functions $a(\lambda)$, 
$b(\lambda)$, and $c(\lambda)$ are in $\FA_r^+$ and satisfy the 
following conditions:
\begin{description}
\item[a)] $a^2$, $b^2$, $c^2$ are rational,
\item[b)] $a$ is odd in $\lambda$, $b$ and $c$ are even in $\lambda$,
\item[c)] $a^2+bc=1$.
\item[d)] $a^2$ is real on $S^1$ and $0\leq a^2\leq 1$ on $S^1$,
\item[e)] $c^2=(b^2)^\ast$.
\end{description}
Furthermore, there exists an odd function $f_+$ in $\FA_r^+$, s.t.\
with $2\halpha(\lambda)=\hp(\lambda)+\hp(-\lambda)$,
$2\hbeta(\lambda)=\hp(\lambda)-\hp(-\lambda)$,
\BEQ \label{pfplus}
\hp=(1+\epsilon_kr_1(\lambda\inv-\lambda))^{|k|}
(1+\epsilon_lr_2(\lambda\inv+\lambda))^{|l|}e^{f_+},
\EEQ
where $k=\epsilon_k|k|,l=\epsilon_l|l|$, we have
\begin{description}
\item[a')] $\halpha$ and $\hbeta^2$ are rational and holomorphic on $\cstar$,
\item[b')] $\halpha$ and $\hbeta^2$ are real on $S^1$,
\item[c')] $\hbeta^2$ is non-positive on $S^1$,
\item[d')] the functions $\hbeta a$, $\hbeta b$, and $\hbeta c$
are rational and holomorphic on $\cstar$. 
\end{description}
The matrix function
$\chi(\lambda)$ is given by 
$\Delta_+^{|k|}\Delta_-^{|l|}\chi=\halpha I+\hbeta h_+Ah_+\inv$.
}

\Proof
By the results of the last sections we know for the functions 
$a,b,c,d$, defined by $S=h_+Ah_+\inv=\tmatrix abcd$:
\begin{itemize}
\item $d=-a$: \bref{deq},
\item $a^2$, $b^2$, $c^2$ are rational functions:
\bref{a2rational}, \bref{b2c2rational},
\item $a$ is odd, $b$ and $c$ are even in $\lambda$: \bref{abparity},
\item $a^2+bc=1$: \bref{abcrel},
\item $a^2$ is real on $S^1$ and $0\leq a^2\leq 1$ on $S^1$:
\bref{a2rational}, \bref{aonS},
\item $c^2=(b^2)^\ast$: \bref{b2c2rel}.
\end{itemize}
Since $(k,l)\in\Sym(\Psi_{mn})$ we have
$F_{m+k,n+l}(\lambda)=\chi(\lambda)F_{mn}(\lambda)$ by 
Theorem~\ref{framesymmetry}.
Moreover,
\BEQ
\Delta_+^{|k|}\Delta_-^{|l|}\chi(\lambda)=\halpha I+\hbeta h_+Ah_+\inv
\EEQ
by~\bref{hchihahS}. By anti-/symmetrization
of~\bref{deus} we get $\halpha(\lambda)=\hp(\lambda)+\hp(-\lambda)$, and
$\hbeta(\lambda)=\hp(\lambda)-\hp(-\lambda)$, where $\hp$ is given
by~\bref{pfplus}. Furthermore,
\begin{itemize}
\item $\halpha$ and $\hbeta^2$ are rational and holomorphic on
$\cstar$: \bref{halphahbetarational},
\item $\halpha$ and $\hbeta^2$ are real on $S^1$: \bref{halphareal},
\bref{alpha2beta2real},
\item $\hbeta^2$ is non-positive on $S^1$: \bref{b2nonpositive}.
\end{itemize}
Finally, $\ha=\hbeta a$, $\hb=\hbeta b$, and $\hc=\hbeta c$ 
are rational and holomorphic on $\cstar$ by~\bref{hSentriesrational}.
\QED

\newsection{} \label{sufftheorem}
We have seen in Section~\ref{lastnecc} under what conditions on $a,b,c$, and
$d$ we can find an $h_+\in\Lambda_r^+\LieSL(2,\Bcc)_\sigma$, s.t.\
$\tmatrix abcd=h_+Ah_+\inv$. This then defines a CMC-immersion $\Psi$
via dressing of the trivial solution with $h_+$. In this section we
characterize those $a,b,c,d,\halpha,\hbeta$ such that a given
$(k,l)\in\Bii^2\setminus\{(0,0)\}$ is in $\Sym(\Psi_{mn})$.

\mytheorem{} {\em
Let there be given three even rational functions $a^2(\lambda)$,
$b^2(\lambda)$, and $c^2(\lambda)$, which for some $r_1,r_2\in\Rplus$
and $\rmin(r_1,r_2)$ defined in Proposition~\ref{loopgroups}
satisfy the following conditions
\begin{description}
\item[a)] $a^2$ is real on $S^1$ and $0\leq a^2\leq 1$ on $S^1$,
\item[b)] $c^2=(b^2)^\ast$,
\item[c)] There exists an $\rmin<r<1$, s.t.\ the restrictions of 
$a^2$, $b^2$, and $c^2$ to $C_r$ are the squares of
functions $a$, $b$, $c$ in $\FA_r^+$,
\item[d)] $a$ is odd, $b$ and $c$ are even in $\lambda$,
\item[e)] $a^2+bc=1$,
\item[f)] $b$ is the square of a holomorphic function on the closure
of $I^{(\rmin)}$.
\end{description}
In addition, with $r$ as in c), 
we assume that there exists an odd function $f_+$ in
$\FA_{r^\prime}^+$, $\rmin<r^\prime\leq r$, 
such that for $\hp=(1+\epsilon_kr_1(\lambda\inv-\lambda))^{|k|}
(1+\epsilon_lr_2(\lambda\inv+\lambda))^{|l|}e^{f_+}$,
$\halpha(\lambda)=\frac12(\hp(\lambda)+\hp(-\lambda))$,
$\hbeta(\lambda)=\frac12(\hp(\lambda)-\hp(-\lambda))$, we have
\begin{description}
\item[a')] $\halpha$ and $\hbeta^2$ are rational and holomorphic on $\cstar$,
\item[b')] $\halpha$ and $\hbeta^2$ are real on $S^1$,
\item[c')] $\hbeta^2$ is non-positive on $S^1$.
\item[d')] The functions $\hbeta a$, $\hbeta b$, and $\hbeta c$ are
rational and holomorphic on $\cstar$.
\end{description}
Then there exists $0<\rmin<r^{\dprime}\leq r^\prime$ and
$h_+\in\Lambda_{r^\dprime}^+\LieSL(2,\Bcc)_\sigma$, such that
$h_+Ah_+\inv=\tmatrix abc{-a}$. Moreover, for the extended frame
$F_{mn}(\lambda)$ defined by
$h_+(\Unull)^m(\Vnull)^n=F_{mn}(\lambda)p_+(m,n,\lambda)$,
$|\lambda|=r^\dprime$,
we have $F_{m+k,n+l}(\lambda)=\chi(\lambda)F_{mn}(\lambda)$,
where $\chi=\alpha I+\beta h_+Ah_+\inv$ is holomorphic on $\cstar$
and takes values in
$\LieSU(2)$ on $S^1$. In particular, $(k,l)\in\Sym(\Psi_{mn})$ for the
CMC-immersion $\Psi_{mn}$ associated with $F_{mn}(\lambda)$ via Sym's formula.
}

\Proof
Assume, that we have functions $a^2,b^2,c^2$ and $f_+,\hp,\halpha,\hbeta$,
such that a)--f), a')--d') are satisfied. We first want to apply
Theorem~\ref{lastnecc}. We set $d=-a$ and know~\bref{abcinth} by e). 
Since $a,b,c$ are defined at $\lambda=0$ and since $a$ is odd we have
$a(0)=0$, whence $b(0)\neq0$. Thus, by Theorem~\ref{lastnecc}, there
exists some $\rmin(r_1,r_2)<r^\dprime\leq r^\prime<1$ and some
$h_+\in\Lambda_{r^\dprime}^+\LieSL(2,\Bcc)_\sigma$, s.t.\
$S=h_+Ah_+\inv=\tmatrix abc{-a}$. Next we consider the extended
frame defined by the $r$-dressing
$h_+(\Unull)^m(\Vnull)^n=F_{mn}(\lambda)p_+(m,n,\lambda)$ of the
cylinder. Recall that we use in this paper the unique Iwasawa splitting
discussed in Section~\ref{loopgroups}. We also set
$\hchi=\Delta_+^{|k|}\Delta_-^{|l|}\chi=\halpha I+\hbeta S$.
From a') and d') it follows, that $\hchi$ is defined and
holomorphic on $\cstar$. A simple calculation using the parity of
$f_+$ gives
\BEQ
\halpha(\lambda)^2-\hbeta(\lambda)^2=\hp(\lambda)\hp(-\lambda)
=\Delta_+^{|k|}\Delta_-^{|l|}.
\EEQ
Thus, $\det\chi=1$.
Since $\halpha$ is even and $\hbeta$ is odd
in $\lambda$, $\hchi$ and $\chi$ both satisfy the twisting
condition~\bref{DPWtwistcond}, whence
$\chi\in\Lambda_{r^\dprime}\LieSL(2,\Bcc)_\sigma$. As outlined in
Section~\ref{Sdef}, $\chi$ is unitary on $S^1$, iff~\bref{halphareal}
and~\bref{hSreality} are satisfied. But~\bref{halphareal} follows
from a'), b'), and the first part of~\bref{hSreality} is just d'). The
second condition is $(\hbeta a)^\ast=-(\hbeta a)$ and $(\hbeta
c)^\ast=-(\hbeta b)$. To verify this condition we square $\hbeta a$,
$\hbeta b$, and $\hbeta c$ and obtain $((\hbeta a)^2)^\ast=(\hbeta a)^2$
and $((\hbeta c)^2)^\ast=(\hbeta b)^2$, since $\hbeta^2$ and $a^2$ are
real by b'), and $(c^2)^\ast=b^2$ by b). Hence $(\hbeta a)^\ast=\pm\hbeta
a$ and $(\hbeta c)^\ast=\pm\hbeta b$. If $(\hbeta a)^\ast=\hbeta a$, then
$\hbeta a$ is real on $S^1$ and $\hbeta^2a^2=(\hbeta a)^2$ is
non-negative on $S^1$. But a) and c') imply that $\hbeta^2a^2$ is
non-positive on $S^1$. The only possibility for both conditions to hold
is $\hbeta a=0$ on $S^1$. But in this case, of course, also $(\hbeta
a)^\ast=-\hbeta a$ as desired. For the remaining case we consider e)
and obtain $\hbeta^2=\hbeta^2a^2+(\hbeta b)(\hbeta c)$. If $(\hbeta
c)^\ast=+\hbeta b$, then $\hbeta^2=\hbeta^2a^2+|\hbeta c|^2$ on
$S^1$. Hence $|\hbeta c|^2=\hbeta^2(1-a^2)$ implies $\hbeta^2\equiv0$ or
$1-a^2\leq0$. The first case is not possible in view of the form of
$\hp$. The second case yields in view of a), that $a^2\equiv1$ on
$S^1$. Hence $a=\pm1$ on $\Bcc$, a contradiction, since $a(0)=0$. Thus
$(\hbeta c)^\ast=-\hbeta b$ as required.

Finally, we show $(k,l)\in\Sym(\Psi_{mn})$. To this end we multiply
$\chi=\frac{1}{\Delta_+^{|k|}\Delta_-^{|l|}}h_+(\halpha I+\hbeta A)h_+\inv$
from the right with
$h_+(\Unull)^m(\Vnull)^n$. A simple calculation like the one leading
to Lemma~\ref{Ddigression} yields
\BEQ
\frac{1}{\Delta_+^{|k|}\Delta_-^{|l|}}(\halpha I+\hbeta A)
=(\Unull)^k(\Vnull)^lw_+,
\EEQ
where
\BEQ
w_+=\tmatrix{\cosh(f_+)}{\sinh(f_+)}{\sinh(f_+)}{\cosh(f_+)}
\EEQ
is holomorphic on $I^{(r)}$. Thus,
\BEQ
\chi h_+(\Unull)^m(\Vnull)^n=h_+(\Unull)^{m+k}(\Vnull)^{n+l}w_+.
\EEQ
Using the definition of $F_{mn}(\lambda)$ gives
\BEQ
\chi(\lambda)F_{mn}(\lambda)p_+(m,n,\lambda)
=F(m+k,n+l,\lambda)p_+(m+k,n+l,\lambda)w_+.
\EEQ
This shows,
\BEQ
F_{m+k,n+l}(\lambda)=\chi(\lambda)F_{mn}(z,\lambda),
\EEQ
since the Iwasawa splitting chosen is unique and $\chi$ is unitary.
\QED

\section{Hyperelliptic Curves}
\label{algebrogeometric} \message{[algebrogeometric]} 
For fixed lattice constants $r_1,r_2\in\Rplus$, let
$\Psi_{mn}:\Bii^2\rightarrow\threespace$ be a CMC-immersion in the
$r$-dressing orbit of the cylinder.
I.e., if $F_{mn}(\lambda)$ is the
extended frame of $\Psi_{mn}$, then there exists $\rmin(r_1,r_2)<r<1$ and
$h_+\in\Lambda_r^+\LieSL(2,\Bcc)_\sigma$, such that $F_{mn}(\lambda)$ is
given by~\bref{Fdressingdef}. In this chapter we also assume, that
$\Psi_{mn}$ has a periodic metric, i.e., that there exists 
$(k,l)\in\Sym(\Psi_{mn})$, $(k,l)\neq(0,0)$.

As we have seen in the last section, we can characterize discrete
CMC-surfaces with periodic metric in terms of the functions $a,b,c,d$
given by $h_+Ah_+\inv=\tmatrix abc{-a}$. This characterization is in
far reaching analogy to the description of CMC-immersions with
periodic metric in the continuous case~\cite{DoHa:3}.

It is only natural to try to adapt also the algebro-geometric description
of the periodicity conditions to the dicrete case, thereby aiming at
a classification of discrete CMC-tori similar to the one of Pinkall and
Sterling~\cite{PiSt:1}.

First we note, that the properties of the rational functions
$a^2$, $b^2$, $c^2$ stated in Theorem~\ref{nectheorem} a)--e) are the same
as those stated in a)--e) of \cite[Theorem~3.6]{DoHa:3}. It is
therefore natural to introduce a hyperelliptic surface associated to
a discrete periodic CMC-surface in the same way as in the continuous
case.

\newsection{} \label{FCintro}
We define the new variable $\nu=\lambda^2$. We will regard the even
rational functions $a^2(\lambda)$, $b^2(\lambda)$, $c^2(\lambda)$ as
rational functions of $\nu$. Since by b) in Theorem~\ref{nectheorem},
$a$ is an odd function in $\lambda$, $a^2(\lambda)$ has a zero of
order $2(2n-1)$, $n>0$, at $\lambda=0$. Hence, as a function in $\nu$,
$a^2$ has a zero of odd order $2n-1$ at $\nu=0$.
Let $\nu_1,\ldots,\nu_k$ be the points in the $\nu$-plane where $a^2$ has
a pole of odd order, and let $\nu_{k+1},\ldots,\nu_{k+l}$ be the points in the
$\nu$-plane away from $\nu=0$ where $a^2$ has a zero of odd order.

As in the continuous case we easily get

\mylemma{} {\em None of the points $\nu_1,\ldots,\nu_{k+l}\in\cstar$ 
defined above lies on the unit circle.
}

\Proof
We have to show, that $a^2(\lambda)$ has neither a pole nor a zero of
odd order on $S^1$. By d) in
Theorem~\ref{nectheorem}, we know, that $a^2$ has no poles on
$S^1$. By c) and e) in Theorem~\ref{nectheorem}, we have
$(1-a^2)^2=|b^2|^2$ on $S^1$, which shows that $b^2$ and
$c^2=(b^2)^\ast$ are
defined on $S^1$ and also, that $a^2$, $b^2$, and $c^2$ cannot vanish
simultaneously on $S^1$. If $a^2$ has a zero of odd order at
$\lambda_0\in S^1$, then $b^2(\lambda_0)\neq 0$. By d') of
Theorem~\ref{nectheorem}, $(\hbeta a)^2$
and $(\hbeta b)^2$ are squares of holomorphic functions on $\cstar$. Thus, the
function $\hbeta^2=\frac{(\hbeta a)^2}{a^2}=\frac{(\hbeta b)^2}{b^2}$ has
both a zero of odd order and of even order at $\lambda_0$. This
implies $\hbeta\equiv0$, contradicting~\bref{betanontrivial},
\QED

\myprop{} {\em Let $k,l$ and $\nu_1,\ldots,\nu_{k+l}$ be defined as
above. Then $g=\frac12(k+l)$ is an integer and we can order the points
$\nu_1,\ldots,\nu_{2g}$, such that
\BEQ \label{branchpoints}
\nu_{2n}=\tau(\nu_{2n-1}),\kern1cm |\nu_{2n-1}|<1,\kern1cm
n=1,\ldots,g.
\EEQ
}

\Proof
By Lemma~\ref{FCintro}, we have $|\nu_{2n-1}|\neq1$ for
$n=1,\ldots,k+l$. Since $a^2(\nu)$ is real on $S^1$, we have using
Section~\ref{involutions}, that the
set $B=\{\nu_1,\ldots,\nu_{2g}\}$ is invariant under the antiholomorphic
involution $\tau:\nu\rightarrow\nuquer\inv$. 
Since $\tau$ has no fixed points off the unit
circle, we get that $B$ consists of pairs $(\nu_n,\tau(\nu_n))$.
This shows, that $k+l$ is even, whence $g=\frac12(k+l)$ is an
integer, and that we can order $\{\nu_1,\ldots,\nu_{2g}\}$, such
that~\bref{branchpoints} holds.
\QED

\separate
In the following we will order the points $\nu_1,\ldots,\nu_{2g}$
always such that~\bref{branchpoints} holds.

We proceed as in the continuous case by considering the algebraic equation
\BEQ \label{FCdef}
\mu^2=\nu\prod_{k=1}^{2g}(\nu-\nu_k).
\EEQ

\mytheorem{} {\em
The plane affine curve $\tilde{\FC}$
defined by~\bref{FCdef} can be uniquely
extended to a compact Riemann surface $\FC$ of genus $g$. The
meromorphic function $\nu:\tilde{\FC}\rightarrow\Bcc$ extends to a
holomorphic map $\pi:\FC\rightarrow\CPE$ of degree $2$. The
branchpoints of $\pi$ are the roots of $\mu^2$ and the point $\infty$.
}

\Proof
The proof follows immediately from~\cite[Lemma~III.1.7]{Miranda:1}
since $\mu^2$ has odd degree.
\QED

\separate
In other words, \bref{FCdef} is a (nonsingular) hyperelliptic curve,
obtained by compactifying the plane affine curve
$\tilde{\FC}=\FC\setminus\{P_\infty\}$, where
$P_\infty=\pi\inv(\infty)\in\FC$ is a single point.

\myremark{}
In Proposition~\ref{aonFC} it will be shown, that
$a(\nu)=\sqrt{a^2(\nu)}$ can be lifted to a nonconstant meromorphic
function on $\FC$.  It determines the complex structure of $\FC$
uniquely (see \cite[I.1.6]{FaKr:1}, \cite[Satz 8.9]{Fo:1}).

\newsection{} \label{meroonFC}
On $\tilde{\FC}=\FC\setminus\{P_\infty\}$, 
$\nu$ and $\mu$ are holomorphic functions.
Every point on $P\in\FC\setminus\{P_\infty\}$ is determined uniquely
by the values of $\nu$ and $\mu$ at $P$ and can thus be identified
with the pair $(\nu(P),\mu(P))$. Formally, we will also write
$P_\infty=(\infty,\infty)$.
In this notation we get $\pi(\nu,\mu)=\nu$.
Let us define the hyperelliptic involution $I$ on $\FC$ by
\BEQ \label{Idef}
I(\nu,\mu)=(\nu,-\mu).
\EEQ
A point on $\FC$ is a branchpoint iff it is mapped by
$\pi$ to a branchpoint on $\CPE$.
Clearly, the branchpoints of $\FC$
are precisely the fixed points of $I$, i.e., the points $(\nu_k,0)$,
$k=1,\ldots,2g+1$ and $P_\infty$. We will also write $P_0=(0,0)$.
Locally around $P_0$ ($P_\infty$), $\mu$ ($\tmu=\mu\nu^{-(g+1)}$) and
each branch of 
$\lambda(\nu,\mu)=\sqrt{\nu}$ ($\lambda\inv(\nu,\mu)=\nu^{-\frac12}$)
are local coordinates on $\FC$, as can be easily checked using the
representation of $\FC$ given in~\cite[Chapter III]{Miranda:1}.

Let us investigate the set of meromorphic functions on $\FC$.
By~\cite[Proposition~1.10]{Miranda:1}, every meromorphic function on
$\FC$ can be uniquely written as
\BEQ \label{merofuncFC}
f(\nu,\mu)=f_1(\nu)+f_2(\nu)\mu,
\EEQ
with two rational functions $f_1$, $f_2$.

\myremark{} 
It is clear from the representation~\bref{merofuncFC} 
of meromorphic functions on
$\FC$, that each rational function $f_1(\nu)$ can be lifted to a
meromorphic function on $\FC$ by setting
$f(\nu,\mu)=f_1(\nu)$. Clearly, then $f\circ I=f$ for such a function.

Conversely, if $f:\FC\rightarrow\Bcc$ is meromorphic, then it can be
identified with a rational function $f_1(\nu)$ iff $f_2(\nu)\equiv0$
in~\bref{merofuncFC}, i.e., iff it satisfies
$f\circ I=f$. We will frequently 
use this identification of rational functions in
$\nu$ with $I$-invariant meromorphic functions on $\FC$.

\newsection{} \label{FCinvolutions}
Let us define
\BEQ \label{FSdef}
\FS=\pi\inv(S^1)=\{(\nu,\mu)\in\FC;\nu\in S^1\}.
\EEQ
The set $\FS$ is connected if $g$ is even, and has two connected 
components if $g$ is odd.
Since $\FS$ is contained in $\tilde{\FC}$, we can identify it with a
subset of $\Bcc^2$.
Using the antiholomorphic involution 
$\tau:\nu\rightarrow\nuquer\inv$ defined in Section~\ref{involutions},
we define the map
$\tsigma:\tilde{\FC}\setminus\{P_0\}\rightarrow\tilde{\FC}\setminus\{P_0\}$ by
\BEQ \label{tauonFC}
\tsigma:(\nu,\mu)\longmapsto
(\nuquer\inv,\nuquer^{-(g+1)}\left(\prod_{i=1}^{2g}\nu_j\right)^{\frac12}
\overline{\mu}).
\EEQ
We will choose the sign of the square root such that the points on
$\FS$ are fixed by $\tsigma$.

The following is well known (see e.g.~\cite{Ja:1}):

\mytheorem{} {\em
The map $\tsigma$ defined by~\bref{tauonFC} can
be extended to an antiholomorphic involution $\FCinvolution$ on $\FC$,
which preserves the points of $\FS\subset\FC$.

Furthermore, $\FCinvolution$ commutes with the hyperelliptic
involution and leaves invariant the set of branchpoints of $\FC$.
}

\Proof
Using Theorem~\ref{FCintro} it is easily checked, that $\tsigma$
defines an antiholomorphic involution on $\tilde{\FC}\setminus\{P_0\}$.
By using the appropriate branches of $\lambda$ and $\lambda\inv$
in the vicinity of $P_0$ and $P_\infty$ on $\FC$, we get
$\tsigma(\lambda)=\lambdaquer\inv$ in local coordinates around $P_0$
and $P_\infty$, whence $\tsigma$ extends to an antiholomorphic involution
$\FCinvolution$ on $\FC$, which maps $P_0$ to $P_\infty$. By the
choice of the square root in~\bref{tauonFC}, $\FCinvolution$ fixes the
points on $\FS$.
$\FCinvolution$ clearly commutes with $I$. If $P$ is a branchpoint of
$\FC$, then $I(P)=P$. Therefore,
$I(\FCinvolution(P))=\FCinvolution(I(P))=\FCinvolution(P)$ and
$\FCinvolution(P)$ is also a branchpoint. 
Thus, $\FCinvolution$ leaves invariant the set of branchpoints of $\FC$.
\QED

\separate
For a scalar function on $\FC$ we also define
\BEQ \label{FCstardef}
f^\ast=\overline{f\circ\FCinvolution}.
\EEQ
Since, by Proposition~\ref{FCinvolutions},
$\FCinvolution$ fixes the points of $\FS$, we get

\mylemma{} {\em
Let $f$ be a meromorphic function defined on a
$\FCinvolution$-invariant subset of $\FC$ which contains $\FS$. 
Then $f$ is real on $\FS$ iff $f^\ast=f$ and
$f$ is purely imaginary on $\FS$ iff $f^\ast=-f$.
}

\newsection{} \label{aonFC}
Let us now investigate the properties of $a^2$ w.r.t.\ $\FC$.

\myprop{} {\em
The rational function $a^2$ defined in Theorem~\ref{nectheorem} is of
the form
\BEQ
a^2(\nu)=f(\nu)^2\mu^2(\nu),
\EEQ
where $f$ is rational and defined at $\nu=0$.
The function $a=f\mu$ is a meromorphic function on 
$\FC$, which satisfies
\BEQ \label{FCaodd}
a\circ I=-a
\EEQ 
and 
\BEQ \label{FCareal}
a^\ast=a.
\EEQ
}

\Proof
By the definition of $\mu^2$, the quotient $\frac{a^2}{\mu^2}$ is
rational and has only poles and zeroes of even order. Therefore, it is
the square of a rational function $f(\nu)$. Since $a^2$ has a zero of
odd order at $\nu=0$, $f(\nu)$ is defined at $\nu=0$.
By~\bref{merofuncFC}, $a=f\mu$ is a meromorphic function on $\FC$.
Since $a^2$ is real and non-negative on $S^1$, the function
$a=\sqrt{a^2}$ takes
real values over $S^1$ on $\FC$, whence, by Lemma~\ref{FCinvolutions},
\bref{FCareal} holds.
Furthermore, $a\circ I=-f\mu=-a$ and~\bref{FCaodd} holds.
\QED

\mytheorem{} {\em
Let $\phi$ be a meromorphic function on $\FC$. Then $\phi$ is
antisymmetric w.r.t.\ the hyperelliptic involution 
$I$, i.e.\ $\phi\circ I=-\phi$, iff $\phi(\nu,\mu)=f(\nu)\mu$, where
$f$ is a rational function. Furthermore, in this case the following holds:
\begin{description}
\item[1.] Locally around $P_0$, $\phi$ is an odd
meromorphic function of the local coordinate $\lambda$.
\item[2.] The product $\phi a$ can be identified with a rational function of
$\nu$ on $\CPE$.
\end{description}
}

\Proof
The equivalence statement follows immediately from the
representation~\bref{merofuncFC} of
meromorphic functions on $\FC$. Since $\lambda$
is a local coordinates around $P_0$,
and since $\lambda\circ I=-\lambda$, 1.\ holds.
Finally, by~\bref{FCaodd},
we have that $\phi a$ is invariant under $I$, whence, by
Remark~\ref{meroonFC}, can be identified
with a meromorphic function $(\phi a)(\nu)$ on $\CPE$.
\QED

\newsection{} \label{FCprime}
Let us also define the Riemann surface $\FC^\prime$ on which
$\sqrt{a^2(\lambda)}$ is meromorphic. To be precise (see~\cite[Lemma
III.1.7]{Miranda:1}), let $\FC^\prime$
be the hyperelliptic curve associated to the algebraic equation
\BEQ \label{FCprimedef}
\tmu^2=\prod_{i=1}^{2g}(\lambda-\sqrt{\nu_i})(\lambda+\sqrt{\nu_i}).
\EEQ
The holomorphic map $\lambda$ extends to a holomorphic map
$\pi^\prime:\FC^\prime\rightarrow\Bcc$ of degree $2$.

As for $\FC$ we will identify the points of $\FC^\prime$ with pairs
$(\lambda,\tmu)$.
It should be noted, that $\FC^\prime$ has no branchpoints over
$\lambda=0$ and $\lambda=\infty$ since $a^2$ is even in $\lambda$. In
particular, $(\pi^\prime)\inv(\infty)$ consists of two different
points $P_\infty^{(1)}$ and $P_\infty^{(2)}$. 
However, this will not cause any confusion here.
By $P_0^{(1)}$ and $P_0^{(2)}$ we will denote the two covering points
of $\lambda=0$. Clearly, $\lambda$ is a local coordinate around
$P_0^{(1)}$ and $P_0^{(2)}$ and $\lambda\inv$ is a local coordinate
around $P_\infty^{(1)}$ and $P_\infty^{(2)}$.

Every meromorphic function $\tf$ on $\FC^\prime$ is of the form
\BEQ \label{merofuncFCprime}
\tf(\lambda,\tmu)=\tf_1(\lambda)+\tf_2(\lambda)\tmu
\EEQ
with two rational functions $\tf_1$ and $\tf_2$.

\myremark{} Let $I^\prime$ be the hyperelliptic involution on
$\FC^\prime$. As in Remark~\ref{meroonFC}, we will use the
representation~\bref{merofuncFCprime} to identify rational functions
of $\lambda$ with $I^\prime$-invariant meromorphic functions on 
$\FC^\prime$. In addition, we have

\myprop{} {\em
There is a one-to-one correspondence between
\begin{description}
\item[1.] meromorphic functions on $\FC$ and
\item[2.] meromorphic functions $\tf$ on $\FC^\prime$ which satisfy
\BEQ \label{FCptoFC}
\tf(-\lambda,\tmu)=\tf(\lambda,-\tmu).
\EEQ
\end{description}
}

\Proof
In terms of $\lambda$ we have $\mu^2(\lambda)=\lambda\tmu^2(\lambda)$.
We can therefore rewrite~\bref{merofuncFC} as
\BEQ
f(\nu,\mu)=f_1(\nu)+f_2(\nu)\lambda\tmu
=\tf_1(\lambda)+\tf_2(\lambda)\tmu=\tf(\lambda,\mu),
\EEQ
where $\tf_1(\lambda)=f_1(\nu=\lambda^2)$ and
$\tf_2(\lambda)=f_2(\nu=\lambda^2)\lambda$ are rational. This shows,
that $f$ can be identified with the meromorphic function $\tf$ on
$\FC^\prime$. Since $\tf_1$ is even and
$\tf_2$ is odd in $\lambda$, also~\bref{FCptoFC} is satisfied.
Conversely, if $\tf$ is a meromorphic function on $\FC^\prime$ such
that~\bref{FCptoFC} is satisfied, then in~\bref{merofuncFCprime},
$\tf_1$ is even and $\tf_2$ is
odd in $\lambda$. Therefore, $f_1(\nu)=\tf_1(\lambda^2)$ and
$f_2(\nu)=\lambda\inv\tf_2(\lambda)$ are rational and define a
meromorphic function $f$ on $\FC$ via~\bref{merofuncFC}.
\QED

\mycorollary{} {\em
The function $a^2(\lambda)$ is the square of a meromorphic function on
$\FC^\prime$.
}

\Proof
By Proposition~\ref{aonFC}, $a^2$ is the square of a meromorphic
function $a$ on $\FC$ which is of the form $a=f\mu$, where $f$ is
rational in $\nu$. Using the proof of Proposition~\ref{FCprime} we
can identify $a$ with
a meromorphic function $\ta$ on $\FC^\prime$ which is of the form
$\ta=\tf\tmu$, where $\tf(\lambda)=f(\nu=\lambda^2)\lambda$ 
is an odd rational function in $\lambda$.
Then $\ta^2$ is $I^\prime$-invariant and 
$\ta^2=\tf^2\tmu^2=f^2\mu^2=a^2$, proving the claim.
\QED

\newsection{} \label{betaonFC}
Let us define the non-compact Riemann surface
$\FCstar=\FC\setminus\{P_0,P_\infty\}$.
We already know, that $a$ is meromorphic on $\FC^\prime$ and $\FC$,
and therefore also on $\FCstar$. Now we prove the following important result:

\mytheorem{} {\em 
The functions $\halpha$ and $\hbeta$ are meromorphic on $\FC$ without
poles on $\FCstar$.
The functions $b=\sqrt{b^2}$ and $c=\sqrt{c^2}$ 
are meromorphic on $\FC^\prime$ without poles over $0$ and $\infty$.
}

\Proof
We know by a') and d') in Theorem~\ref{nectheorem}, 
that $\hbeta a$, $\halpha$, and $\hbeta^2$ are even rational functions in
$\lambda$, which don't have poles on
$\cstar$. Therefore, by Remark~\ref{meroonFC}, they can be identified
with $I$-invariant meromorphic functions on $\FC$ without poles on
$\FCstar$. Since also
$a$ is meromorphic on $\FC$, we have that $\hbeta=\frac{\hbeta a}{a}$ is
a meromorphic function on $\FC$. Since the square $\hbeta^2$ is
holomorphic on $\cstar$, $\hbeta$ has no poles on $\FCstar$.

By Proposition~\ref{FCprime}, $\hbeta$ can be identified with a
meromorphic function on $\FC^\prime$.
Furthermore, $\hbeta b$ and $\hbeta c$ are rational functions of
$\lambda$. Thus, $\hbeta b$, $\hbeta c$, $b=\frac{\hbeta b}{\hbeta}$ and 
$c=\frac{\hbeta c}{\hbeta}$ are meromorphic on $\FC^\prime$.
Since, by Theorem~\ref{nectheorem}, $b$ and $c$ are in $\FA_r^+$, they
can be continued holomorphically to $\lambda=0$ on $\CPE$. By e) in
Theorem~\ref{nectheorem}, the same holds for $b$ and $c$ around
$\lambda=\infty$.  Since $\lambda$ is a local coordinate around
$P_0^{(1)}$ and $P_0^{(2)}$ and $\lambda\inv$ is a local coordinate
around $P_\infty^{(1)}$ and $P_\infty^{(2)}$, the functions $b$ and
$c$ can be extended holomorphically to these points, which finishes
the proof.
\QED

The proof of the following proposition is the same as the one of
\cite[Prop.~4.6]{DoHa:3}.

\myprop{} {\em
With $(\cdot)^\ast$ defined in Section~\ref{FCinvolutions}, the holomorphic
functions $\halpha$ and $\hbeta$ on $\FCstar$ satisfy
\BEQ
\halpha^\ast=\halpha,\kern3cm\hbeta^\ast=-\hbeta,
\EEQ
\BEQ
\mbox{$\alpha\circ I=\alpha$ and $\beta\circ I=-\beta$.}
\EEQ
Furthermore, the meromorphic functions $\alpha$, $\beta$ and $b$ and
$c$ on $\FC^\prime$ satisfy
\BEQ
c=b^\ast,
\EEQ
\BEQ
\mbox{$\alpha\circ I^\prime=\alpha$ and $\beta\circ I^\prime=-\beta$,}
\EEQ
\BEQ
\mbox{$b\circ I^\prime=-b$ and $c\circ I^\prime=-c$.}
\EEQ
}
%

It is easily possible to pursue this road of deriving results analogous to the
continuous case, as it was done in
\cite[Sections~4.8--4.9]{DoHa:3}. In particular, one can define higher
commuting flows acting on discrete CMC-surfaces and finite type
solutions of the discrete $\sinh$-Gordon equation in 
precisely the same way as
in the continuous case. The same argument as in the proof of
\cite[Theorem~4.9]{DoHa:3} then shows, that discrete CMC-surfaces with
periodic metric belong to finite type solutions of the discrete
$\sinh$-Gordon equation. However, since there is up to now no theory
of finite type solutions of the discrete $\sinh$-Gordon equation, this
does not promise any help towards our goal.

\section{Algebro-geometric description of surfaces with periodic
metric} \label{tori} \message{[tori]}
For a discrete CMC-surface with periodic metric we defined in 
Section~\ref{FCintro} a nonsingular hyperelliptic curve $\FC$.
In this section, we will show, that $\FC$ allows us to express the periodicity
conditions for discrete CMC-surfaces stated in Theorem~\ref{nectheorem} and
Theorem~\ref{sufftheorem} in terms of algebro-geometric data.

We will also investigate the case, that the discrete CMC-surface
$\Psi_{mn}$ under consideration does not only have a periodic or doubly
periodic metric, but closes in $\threespace$, i.e., that the image
$\{\Psi_{mn};(m,n)\in\Bii^2\}\subset\threespace$ consists of finitely
many points. This case is an obvious analogue of a CMC-torus. We will
call such surfaces {\em discrete CMC-tori}.

\newsection{} \label{cycles}
We will first reformulate the statement of Theorem~\ref{nectheorem} in
terms of algebro-geometric data. We start with the same
assumptions as in Section~\ref{algebrogeometric}: For fixed lattice
constant $r_1,r_2\in\Rplus$, let
$\Psi_{mn}:\Bii^2\rightarrow\threespace$ be a discrete CMC-surface
with extended frame $F_{mn}\in\FF_0(r_1,r_2)$, obtained by dressing
the discrete cylinder with some
$h_+\in\Lambda_r^+\LieSL(2,\Bcc)_\sigma$, $\rmin(r_1,r_2)<r<1$. We also
assume, that
$\Sym(\Psi_{mn})$ contains a nontrivial element $(k,l)\neq(0,0)$.
I.e., $\Psi_{mn}$ has periodic metric and we can define the
hyperelliptic curve $\FC$ as in Section~\ref{FCintro}.

We introduce a standard homotopy
basis for $\FC$ which is adapted to the $\FCinvolution$-symmetry of
$\FC$ stated in Proposition~\ref{FCinvolutions}.
Let $a_1,\ldots,a_g,b_1,\ldots,b_g$, $g$ the genus
of $\FC$, be a canonical basis of $H_1(\FC,\Bii)$, such that the
intersection numbers are given by
\BEQ
a_ia_j=0,\kern1cm
b_ib_j=0,\kern1cm a_ib_j=\delta_{ij},\kern1cm i,j=1,\ldots,g.
\EEQ
For the cycles $a_k$ we choose (see~\cite[VII.7.1]{FaKr:1})
\BEQ
a_k=\gamma_k-I\circ\gamma_k,
\EEQ
where $\gamma_k$ is a curve joining 
the branchpoints over $\nu_{2k-1}$ and $\nu_{2k}$, which satisfies
$\FCinvolution\circ\gamma=-\gamma$. Then
\BEQ \label{anachI}
I\circ a_k=-a_k
\EEQ
and
\BEQ \label{aktau}
\FCinvolution\circ a_k=-a_k,
\EEQ
since $\FCinvolution$ and $I$ commute.
I.e., the cycles $a_k$ are up to orientation invariant under
$\FCinvolution$ and $I$. In addition, we can choose $b_k$ such that
\BEQ \label{bktau}
\FCinvolution\circ b_k=b_k-a_k+\sum_{j=1}^ga_j.
\EEQ

\newsection{}
By Theorem~\ref{betaonFC}, $\halpha$ and $\hbeta$ are meromorphic
functions on $\FC$ without poles on $\FCstar$. Therefore,
$\diff\halpha$ and $\diff\hbeta$ and $\halpha\diff\hbeta-\hbeta\diff\halpha$
are meromorphic one-forms on $\FC$ without poles on $\FCstar$.
We define the meromorphic one-form $\omega$ on $\FC$ by
\BEQ \label{dpdef}
\omega:=\frac{\halpha\diff\hbeta-\hbeta\diff\halpha}{
\Delta_+^{|k|}\Delta_-^{|l|}}.
\EEQ
Using~\bref{hahbdef} and~\bref{a2b2}, this can be rewritten as
\BEQ \label{dpln}
\omega=\alpha\diff\beta-\beta\diff\alpha=(\alpha-\beta)\diff(\alpha+\beta)
=\frac{\diff(\alpha+\beta)}{\alpha+\beta}
\EEQ
for
\BEQ
\alpha=\frac{1}{\Delta_+^{|k|}\Delta_-^{|l|}}\halpha,\kern3cm
\beta=\frac{1}{\Delta_+^{|k|}\Delta_-^{|l|}}\hbeta.
\EEQ
{\em From now on, we restrict our attention to the case that $k$ and $l$
are even.} Since $\Sym(\Psi_{mn})$ is an additive group, it is clear,
that for every discrete CMC-surface with periodic metric, there exist
nontrivial $(k,l)\in\Sym(\Psi_{mn})$ with $k$ and $l$ even.
Thus our choice of $k$ and $l$ even does not impose a restriction on
the class of discrete CMC-surfaces under consideration.
Since $\Delta_+^2$ and $\Delta_-^2$ are even rational functions of
$\lambda$, i.e.~by Remark~\ref{meroonFC}, 
meromorphic functions on $\FC$, it is clear that in this case
with $\halpha$ and $\hbeta$ also $\alpha$ and $\beta$ are
meromorphic on $\FC$.

In the following we will denote by $\FI^{(r)}$ the set of all points
on $\FC$ which are mapped to $I^{(r)}$ by the projection $\pi$.
Since the surface $\FC$ constructed in Section~\ref{FCintro} has no
branchpoints in $\FI^{(r)}\setminus\{P_0\}$, the local coordinate
$\lambda$ extends to all of $\FI^{(r)}$, i.e., $\FI^{(r)}$ is a chart
domain on $\FC$ around $P_0$.
Similarly, we define $\FE^{(\frac{1}{r})}$ to be the set of all points on $\FC$
which are mapped to $\tau(I^{(r)})=\{\nu;|\nu|>\frac{1}{r}\}$.
Obviously, $\FCinvolution(\FI^{(r)})=\FE^{(\frac{1}{r})}$.
Each branch of $\lambda\inv$ extends to a local coordinate on
$\FE^{(\frac{1}{r})}$. For a chosen branch of $\lambda$ on $\FI^{(r)}$
we fix a branch of $\lambda\inv$ by requiring
\BEQ
\mbox{\rm $\lambda\inv(\FCinvolution(P))=\overline{\lambda(P)}$ for
all $P\in\FI^{(r)}$.}
\EEQ
Thus $\FCinvolution$ restricts to the map
$\lambda\mapsto\lambdaquer\inv$ from $\FI^{(r)}$ to $\FE^{(\frac1r)}$.

From~\bref{deusorigin} and~\bref{hpdef}, we get on $\FI^{(r)}$:
\BEQ \label{abinlocallambda}
\alpha+\beta=\frac{(1+\epsilon_k r_1(\lambda\inv-\lambda))^{|k|}
(1+i\epsilon_lr_2(\lambda\inv-\lambda))^{|l|}}{
(1-r_1^2(\lambda\inv-\lambda)^2)^{\frac{|k|}{2}}
(1+r_2^2(\lambda\inv+\lambda)^2)^{\frac{|l|}{2}}}e^{f_+}
=S_{r_1}(\lambda)^{\frac{k}{2}}T_{r_2}(\lambda)^{\frac{l}{2}}e^{f_+},
\EEQ
where $S_{r_1}$ and $T_{r_2}$ were defined in Section~\ref{chiform}
and $f_+$ is a holomorphic function on $\FI^{(r)}$.
This shows, that $\alpha+\beta$ has no pole at $P_0$. Therefore,
by~\bref{dpln}, $\omega$ has poles on $\FI^{(r)}$ only where $S_{r_1}$
and $T_{r_2}$ have zeroes or poles, 
and all of these poles are simple. Furthermore, by
Proposition~\ref{cylinder} and c') in Theorem~\ref{nectheorem},
\BEQ
(\alpha+\beta)^\ast=\alpha-\beta=(\alpha+\beta)\inv.
\EEQ
Therefore, $\alpha+\beta$ and $\omega$ are 
meromorphic on $\FE^{(\frac1r)}$.

Recall the definition of $\lambda_+$ and $\lambda_-$ in
Lemma~\ref{cylinder}.  On $I^{(r)}$, $S_{r_1}$ has precisely one zero
at $\lambda_+$ and one pole at $-\lambda_+$, both are simple.  Also
on $I^{(r)}$, $T_{r_2}$ has precisely one (simple) zero at
$i\lambda_-$ and one (simple) pole at $-i\lambda_-$.  Therefore, the
representation~\bref{dpln} together with~\bref{abinlocallambda} shows,
that $\omega$ has precisely $4$ poles on $\FI^{(r)}$, which are
located at $\pm\lambda_+$ and $\pm i\lambda_-$ with $\lambda_\pm$
defined as in Lemma~\ref{cylinder}.  All of these poles are simple and
the residues of $\omega$ are given by
\BEA \label{firstset}
\res_{\lambda_+}\omega & = & -\res_{-\lambda_+}\omega=\frac{k}{2},\nonumber\\
\res_{i\lambda_-}\omega & = & -\res_{-i\lambda_-}\omega=\frac{l}{2}.
\EEA
Using the local coordinate $\lambda\inv$ on $\FE^{(\frac1r)}$, we get
another set of $4$ simple poles of $\omega$, which are located at
$\pm\lambda_+\inv$ and $\pm i\lambda_-\inv$. The residues of $\omega$
at these poles are given by
\BEA \label{secondset}
\res_{-\lambda_+\inv}\omega & = & -\res_{\lambda_+\inv}\omega
=\frac{k}{2},\nonumber\\ 
\res_{-i\lambda_-\inv}\omega & = & -\res_{i\lambda_-\inv}\omega
=\frac{l}{2}.
\EEA
Since $\halpha\diff\hbeta-\hbeta\diff\halpha$ 
has no poles on $\FCstar$, the form $\omega$ has by~\bref{dpdef} no
further poles on $\FC$. We have proved the first part of the

\mylemma{} {\em
The one-form $\omega$ has precisely $8$ simple poles on $\FC$ with residues
given by~\bref{firstset} and~\bref{secondset}.
Furthermore,
\BEQ \label{dpimaginary}
\FCinvolution^\ast\omega=-\overline{\omega}
\EEQ
and
\BEQ \label{dpak}
\int_{a_k}\omega=0,\kern1cm k=1,\ldots,g.
\EEQ
$\omega$ is an Abelian differential of the third kind, which is
completely determined by~\bref{firstset}, \bref{secondset}
and~\bref{dpak}.
}

\Proof
Note that by Propositions~\ref{cylinder},\ref{betaonFC} and
Lemma~\ref{FCinvolutions}, we have
\BEQ
\alpha^\ast=\overline{\alpha\circ\FCinvolution}=\alpha,\kern2cm
\beta^\ast=\overline{\beta\circ\FCinvolution}=-\beta.
\EEQ
From this it follows,
\BEQ
\FCinvolution^\ast\diff\alpha=\overline{\diff\alpha},\kern2cm
\FCinvolution^\ast\diff\beta=-\overline{\diff\beta}.
\EEQ
With this, Eq.~\bref{dpimaginary} follows from~\bref{dpdef}.
Since an Abelian differential of the third kind with only simple poles
is completely determined by its residues and its $a_k$-cycles
(see~\cite[Prop.~III.3.3]{FaKr:1}),
it only remains to be shown, that~\bref{dpak} holds.

Since $\omega$ is the logarithmic derivative of the meromorphic
function $\alpha+\beta$ on $\FC$, its integral over a closed cycle on
$\FC$ (avoiding the poles of $\alpha+\beta$) 
is an integer multiple of $2\pi i$. By~\bref{dpimaginary}
and~\bref{aktau}, we have
\BEQ
\overline{\int_{a_k}\omega}=\int_{a_k}\overline{\omega}=
-\int_{a_k}\FCinvolution^\ast\omega
=-\int_{\FCinvolution\circ a_k}\omega=\int_{a_k}\omega.
\EEQ
Thus, $\int_{a_k}\omega$ is real, which implies~\bref{dpak}.
\QED

\newsection{} \label{dOplus}
Let $\Omega_+$ be an Abelian differential of the third kind on $\FC$,
which has two simple poles at $\lambda_+$ and $-\lambda_+$ in
$\FI^{(r)}$ with residue
\BEQ \label{dOpres}
\res_{\lambda_+}\Omega_+=-\res_{-\lambda_+}\Omega_+=1.
\EEQ
Such a differential exists by~\cite[Theorem~II.5.2]{FaKr:1}. If we require
\BEQ \label{dOpak}
\int_{a_k}\Omega_+=0,\kern1cm k=1,\ldots,g,
\EEQ
then $\Omega_+$ is uniquely determined by~\bref{dOpres} and~\bref{dOpak}
(see \cite[Prop.~III.3.3]{FaKr:1}).
We also introduce a second Abelian differential $\Omega_-$ of the third kind
on $\FC$, which has two simple poles at $i\lambda_-$ and $-i\lambda_-$
in $\FI^{(r)}$. It is uniquely determined by
\BEQ \label{dOmres}
\res_{i\lambda_-}\Omega_-=-\res_{-i\lambda_-}\Omega_-=1
\EEQ
and
\BEQ \label{dOmak}
\int_{a_k}\Omega_-=0,\kern1cm k=1,\ldots,g.
\EEQ
If we define
\BEQ
\Omega_+^\ast=\overline{\FCinvolution^\ast\Omega_+},\kern3cm
\Omega_-^\ast=\overline{\FCinvolution^\ast\Omega_-},
\EEQ
then $\Omega_+^\ast$ and $\Omega_-^\ast$ are also Abelian differentials of the
third kind with simple poles at $\pm\lambda_+\inv$ and
$\pm i\lambda_-\inv$ in $\FE^{(\frac1r)}$, respectively.
The residues of $\Omega_+^\ast$ and $\Omega_-^\ast$ are given by
\BEQ \label{dOpstarres}
\res_{\lambda_+\inv}\Omega_+^\ast=-\res_{-\lambda_+\inv}\Omega_+^\ast=1
\EEQ
and
\BEQ \label{dOmstarres}
\res_{i\lambda_-\inv}\Omega_-^\ast=-\res_{-i\lambda_-\inv}\Omega_-^\ast=1.
\EEQ
By~\bref{aktau} and~\bref{dOpak} we get for $k=1,\ldots,g$:
\BEQ \label{dOpstarak}
\int_{a_k}\Omega_+^\ast=\int_{a_k}\overline{\FCinvolution^\ast\Omega_+}
=\overline{\int_{\FCinvolution\circ a_k}\Omega_+}=-\overline{\int_{a_k}\Omega_+}=0.
\EEQ
In the same way it follows from~\bref{dOmak} that
\BEQ \label{dOmstarak}
\int_{a_k}\Omega_-^\ast=0,\kern1cm k=1,\ldots,g.
\EEQ
Comparing~\bref{dpak} with~\bref{dOpak}, \bref{dOmak}, \bref{dOpstarak},
\bref{dOmstarak} and comparing~\bref{firstset}, \bref{secondset}
with~\bref{dOpres}, \bref{dOmres}, \bref{dOpstarres}, \bref{dOmstarres},
we see that
\BEQ \label{dpdO}
\omega=\frac{k}{2}(\Omega_+-\Omega_+^\ast)+\frac{l}{2}(\Omega_--\Omega_-^\ast).
\EEQ
Let us set
\BEQ
U^+_n=\int_{b_n}\Omega_+,\kern3cm U^-_n=\int_{b_n}\Omega_-,\kern1cm n=1,\ldots,g.
\EEQ
Then with~\bref{bktau} and~\bref{dOpak},\bref{dOmak},
\BEQ
\int_{b_n}\Omega_+^\ast=\overline{U^+_n},\kern3cm
\int_{b_n}\Omega_-^\ast=\overline{U^-_n},\kern1cm k=1,\ldots,g.
\EEQ

\mylemma{} {\em
The integrals $U^+_n$ and $U^-_n$ of $\Omega_+$ and $\Omega_-$ over the
$b_n$-cycles satisfy
\BEQ \label{UUcond}
k\Imag(U^+_n)+l\Imag(U^-_n)=2\pi m_n,\kern1cm m_n\in\Bii.
\EEQ
}

\Proof
We can deform the cycles $b_n$ homotopically such that none of them
crosses one of the poles of $\omega$. Then the integral of $\omega$
over $b_n$ is an integer multiple of $2\pi i$, since it is the
logarithmic derivative of the meromorphic function $\alpha+\beta$ on
$\FC$. From this and~\bref{dpdO} the claim follows.
\QED

The question naturally arises, if there really exist algebro-geometric
data $\FC$, s.t.\ the periodicity conditions are satisfied.
In the continuous case this question was answered in the affirmative
first by Ercolani, Kn\"orrer, and Trubowitz for even genus $g>0$ of
the spectral curve, later by Jaggy~\cite{Ja:1} for every genus $g>1$.
It should in principle be possible to start a similar investigation
for the discrete case.

There exists numerical evidence, that as in the continuous case, the
conditions cannot be satisfied for genus $g=1$ of the spectral
curve. However, T.~Hoffmann has constructed discrete Delaunay type
surfaces using a different approach. It remains to be investigated, if
these surfaces are covered by the method of this paper.

Finally, besides of existence proofs, it is an interesting problem to
construct such surfaces numerically, as it was done in the continuous
case by M. Heil~\cite{Heil:1}.
Using theta functions it is no problem to give an explicit formula
for the discrete `immersion' $\Psi$. The formulae look similar to
those of Bobenko, just replacing the differentials of the second kind
by differentials of the third kind.

\end{document}